\documentclass[11pt]{article}
\usepackage{amscd}
\usepackage{amsfonts}
\usepackage{amsmath}
\usepackage{amssymb}
\usepackage{amsthm}
\usepackage{bbm}
\usepackage{CJK}
\usepackage{fancyhdr}
\usepackage{graphicx}
\usepackage{hyperref}
\usepackage{indentfirst}
\usepackage{latexsym}
\usepackage{mathrsfs}
\usepackage{xypic}
\usepackage{amsmath,amssymb,amsfonts,amsthm,fancyhdr}
\usepackage{epsfig,graphicx,picins,picinpar,subfigure}
\usepackage{pstricks}
\usepackage{amssymb}
\usepackage{xypic}
\usepackage[compress]{cite}

\newtheorem{theorem}{Theorem}[section]
\newtheorem{prop}[theorem]{Proposition}
\newtheorem{lemma}[theorem]{Lemma}
\newtheorem{coro}[theorem]{Corollary}
\newtheorem{prop-def}{Proposition-Definition}[section]

\theoremstyle{definition}
\newtheorem{definition}[theorem]{Definition}

\newtheorem{remark}[theorem]{Remark}
\newtheorem{exam}[theorem]{Example}

\usepackage[top=1in,bottom=1in,left=1.25in,right=1.25in]{geometry}
\def\<{\langle}
\def\>{\rangle}

\date{\today}
\begin{document}
\renewcommand{\baselinestretch}{1.2}
\renewcommand{\arraystretch}{1.0}
\title{\bf
Cohomology and deformations of nonabelian  embedding tensors between Lie triple systems }
\author{{\bf     Wen Teng\footnote
        { Corresponding author:~~tengwen@mail.gufe.edu.cn}  }\\
{\small  School of Mathematics and Statistics, Guizhou University of Finance and Economics} \\
{\small  Guiyang  550025, P. R. of China}}
 \maketitle
\begin{center}
\begin{minipage}{13.cm}

{\bf Abstract}
In this paper,  first  we introduce the notion of   nonabelian  embedding tensors between  Lie triple systems and show that nonabelian embedding tensors  induce naturally  3-Leibniz algebras. Next, we construct  an $L_{\infty}$-algebra  whose Maurer-Cartan elements are  nonabelian  embedding tensors. Then, we have the  twisted
$L_{\infty}$-algebra that governs deformations of  nonabelian  embedding tensors.  Following this,  we establish the cohomology of a  nonabelian  embedding tensor  between  Lie triple systems  and realize it as the cohomology of the descendent 3-Leibniz algebra with coefficients in a suitable
representation.  As applications, we consider infinitesimal deformations of a nonabelian embedding tensor between Lie triple systems and demonstrate that they are governed by the above-established cohomology.
Furthermore, the notion
of Nijenhuis elements associated with a nonabelian embedding tensor is
introduced to characterize trivial infinitesimal deformations.
Finally, we provide relationships between nonabelian  embedding tensors
on  Lie algebras and associated Lie triple systems.

 \smallskip

 \smallskip

 \end{minipage}
 \end{center}
 \normalsize\vskip0.5cm

\section{Introduction}
\def\theequation{\arabic{section}. \arabic{equation}}
\setcounter{equation} {0}

The notion  of Lie triple systems  first appeared in Cartan's work \cite{Cartan} on  totally geodesic submanifolds and symmetric spaces,  which were later investigated by  Jacobson from an algebraic
viewpoint in \cite{Jacobson1,Jacobson2}.  Lister constructed a structure theory of Lie triple systems in \cite{Lister}.
After that,   Lie triple system has important applications in physics, such as
quantum mechanics theory  and  numerical analysis of differential equations \cite{Munthe-Kaas}.
 The representation, cohomology and homotopy theory of Lie triple systems were
established in \cite{Hodge,X23,Yamaguti,Kubo}. Lie triple systems  got a lot of attention \cite{Yamaguti1,Hajjaji24,Wu22, Chtioui, Teng23,Teng25}.

The concept of embedding tensors has been extensively studied and applied in numerous areas of mathematics and physics.
Embedding tensors initially emerged in the exploration of gauged supergravity theory \cite{Nicolai}.
  From a   mathematical point of view, embedding tensors are also known as averaging operators \cite{Aguiar}. In particular, an
averaging operator on a Lie algebra gives rise to a Leibniz algebra.  Averaging operators have attracted attention from both mathematics and physics \cite{Miller,Moy,Cao,Das24}.

 Recently, the authors studied  the  controlling $L_\infty$-algebra, cohomology and homotopy of embedding tensors on   Lie algebras in \cite{Sheng}. Embedding tensors on arbitrary types of algebras such as Lie$_\infty$-algebras,  3-Lie algebras, 3-Hom-Lie algebras and Lie triple systems were studied in\cite{Hu,Teng,Caseiro,Teng25}.
As a nonabelian generalization of embedding tensors on Lie algebras, in \cite{Tang} we can find the case where the representation space is also Lie algebra and, in this case, the embedding  tensors are called nonabelian  embedding tensors.
Such generalization for averaging operators on associative algebra was previously addressed in \cite{Pei}, which can be understood as weighted averaging operators.

The purpose of this paper is to introduce the nonabelian embedding tensor   between  Lie triple systems, which is a nonabelian generalization of the embedding tensor on Lie triple systems \cite{Teng25}.
To validate that nonabelian embedding tensors are indeed robust mathematical structures, we use Voronov's higher derived brackets \cite{Voronov} to construct an $L_\infty$-algebra whose Maurer-Cartan
elements are nonabelian embedding tensors. Consequently, we obtain a twisted $L_\infty$-algebra, it is a differential graded Lie algebra governing deformations of a nonabelian embedding tensor.
Furthermore, a nonabelian embedding tensor induces a   3-Leibniz algebra, which is called the descendent 3-Leibniz algebra of
the nonabelian embedding tensor. We introduce a cohomology theory of nonabelian
embedding tensors via the  cohomology of the descendent 3-Leibniz
algebras.   As applications, we use the first   cohomology group to study infinitesimal deformations of a nonabelian
embedding tensor    between   Lie triple systems. In particular, we introduce the notion of a Nijenhuis element associated to a nonabelian embedding tensor, which leads to a trivial infinitesimal deformations of the nonabelian embedding tensor between  Lie triple systems.

The paper is organized as follows.   In Section  \ref{sec: NETs},  we introduce the notion of nonabelian
embedding tensors
between Lie triple systems. We show that nonabelian embedding tensors can be characterized by graphs of the  nonabelian hemisemidirect product 3-Leibniz algebra.
In Section \ref{sec: Maurer-Cartan},  we construct an $L_\infty$-algebra whose Maurer-Cartan elements are nonabelian  embedding tensors
between  Lie triple systems.
In Section \ref{sec: Cohomology},  we introduce a cohomology theory of nonabelian embedding tensors between Lie triple systems.
 In Section \ref{sec: Deformations},  we study  infinitesimal deformations   of a nonabelian embedding tensor  between Lie triple systems
 using the established cohomology theory.
 Finally, in Section \ref{sec: From}, we describe
some connections between nonabelian embedding tensors on Lie algebras and associated Lie
triple systems.

\section{ Nonabelian embedding tensors between Lie triple systems }\label{sec: NETs}
\def\theequation{\arabic{section}.\arabic{equation}}
\setcounter{equation} {0}

In this section,  we introduce the notion of nonabelian embedding tensors between Lie triple systems,  which can be characterized by the graphs of the  nonabelian hemisemidirect product 3-Leibniz algebras.
We work over   a field  $\mathbf{k}$ of characteristic zero.

\begin{definition}  \cite{Jacobson1}
  A Lie triple system   $(\mathbb{L}, [-, -,-])$  is a vector space $\mathbb{L}$ together with  a trilinear operation $[-, -, -]$
on $ \mathbb{L}$    satisfying
\begin{align}
&[x,y,z]+[y,x,z]=0,\label{2.1}\\
&[x,y,z]+[z,x,y]+[y,z,x]=0,\label{2.2}\\
 &[a, b, [x, y, z]]=[[a, b, x],y,z]+ [x,  [a, b, y],z]+ [x,y,[a, b, z]],\label{2.3}
\end{align}
for all $ x, y, z, a, b\in \mathbb{L}$.
 A homomorphism from a  Lie triple system   $(\mathbb{L}', [-, -, -]')$ to a  Lie triple system $(\mathbb{L}, [-, -,-])$ is a linear map $f: \mathbb{L}'\rightarrow \mathbb{L}$ such that
$f([u, v, w]')=[f(u), f(v),f(w)],$ for all $u,v,w\in \mathbb{L}'.$
\end{definition}

\begin{remark}  \cite{Yamaguti} \label{remark: Lie-LTS}
 Let $(\mathbb{L}, [-,-])$ be a     Lie   algebra. Then
$(\mathbb{L}, [-,-,-])$ is a  Lie triple system,
where  $ [x,y,z]=[[x,y],z],$ for   $ x,y,z\in \mathbb{L}.$
Conversely, it was shown in \cite{Yamaguti1} that if $(\mathbb{L}, [-, -,-])$ is a Lie triple system, then
$\Omega(\mathbb{L})=D(\mathbb{L},\mathbb{L})\oplus \mathbb{L}$ with the product
$$[(D(x,y),z),~(D(a,b),c)]=\big(D([a,b,y],x)-D([a,b,x],y)+D(z,c),~D(x,y)c-D(a,b)z\big)$$
is a Lie algebra, where $D(\mathbb{L},\mathbb{L})=\mathrm{span} \{D(x,y)~|~x,y\in \mathbb{L}\}$ and $D(x,y)=[x,y,-]$.
\end{remark}

\begin{definition} \cite{Yamaguti}
A representation of  a Lie triple system  $(\mathbb{L}, [-, -,-])$ on  a vector space $V$ is a bilinear map $\theta: \mathbb{L}\times \mathbb{L}\rightarrow \mathrm{End}(V)$, such that for all  $x,y,a,b\in \mathbb{L}$,
%Let $(\mathbb{L}, [-, -,-])$ be a Lie triple system and $V$ be a vector space. Then $(V; \theta)$
%is called a representation of $(\mathbb{L}, [-, -,-])$ if   the following equations are satisfied (for all  $x,y,a,b\in \mathbb{L}$),
\begin{align}
& \theta(a,b)\theta(x,y)-\theta(y,b)\theta(x,a)-\theta(x, [y,a,b])+D_\theta(y,a)\theta(x,b)=0,\label{2.4}\\
& \theta(a,b)D_\theta(x,y)-D_\theta(x,y)\theta(a,b)+\theta([x,y,a],b)+\theta(a,[x,y,b])=0,\label{2.5}\
\end{align}
 where $D_\theta: \mathbb{L}\times \mathbb{L}$ to $\mathrm{End}(V)$   is given by
\begin{align}
D_\theta(x,y)=\theta(y,x)-\theta(x,y)\label{2.6}.
\end{align}
From Eqs. \eqref{2.5} and \eqref{2.6} we get
\begin{align}
& D_\theta(a,b)D_\theta(x,y)-D_\theta(x,y)D_\theta(a,b)+D_\theta([x,y,a],b)+D_\theta(a,[x,y,b])=0.\label{2.7}
\end{align}
%We denote it by $(V; \theta)$.
\end{definition}

For example, given a Lie triple system $(\mathbb{L}, [-, -,-])$, there is a natural adjoint representation   on  itself.
The corresponding maps $\theta$ and $D_\theta$ are given by
$$\theta(x,y)z:=[z,x,y] ~\text{and}~$$
$$ D_\theta(x,y)z=\theta(y,x)z-\theta(x,y)z=[x,y,z] , ~\text{for}~  x,y,z\in \mathbb{L}.$$

Next, let us give the definition of a nonabelian embedding tensor    from a
Lie triple system to another Lie triple system. We recall a derived algebra of a Lie triple system
$(\mathbb{L}, [-,-,-])$  is a subsystem $ [\mathbb{L},\mathbb{L},\mathbb{L}]$ , and denoted by $\mathbb{L}^1$ The subspace
$$C(\mathbb{L})=\{x\in\mathbb{L}~|~[x,y,z]=0,~\forall y,z\in\mathbb{L}\}$$
is the center of $(\mathbb{L}, [-,-,-])$.

 \begin{definition}
(1) \cite{Wu22}  Let $(\mathbb{L}, [-,-,-])$ and $(\mathbb{L}', [-,-,-]')$ be two Lie triple systems. Let $\theta: \mathbb{L}\times \mathbb{L}\rightarrow \mathrm{End}(\mathbb{L}')$ be a representation of the Lie triple system $\mathbb{L}$ on the vector space $\mathbb{L}'$. If for all $x,y\in\mathbb{L},$
\begin{align}
& \theta(x,y) u\in C(\mathbb{L}'),\label{2.8}\\
& \theta(x,y) [u,v,w]'=0,\label{2.9}
\end{align}
for all $x,y\in \mathbb{L}$ and $u,v,w\in\mathbb{L}'$,
then $\theta$ is called an   action of  $(\mathbb{L}, [-,-,-])$  on  $(\mathbb{L}', [-,-,-]')$.

(2) An action $\theta$ of $ \mathbb{L} $  on  $ \mathbb{L}' $ is called a
coherent action if $\theta$ satisfies
$$\theta(x,y)\in \mathrm{Der}(\mathbb{L}'),~~ \text{for all}~~ x,y\in\mathbb{L}.$$

\end{definition}
Let $\theta$ be a coherent action of $ \mathbb{L} $  on  $ \mathbb{L}' $.  By  Eqs. \eqref{2.6}, \eqref{2.8} and \eqref{2.9}, we  deduce that
\begin{align}
& D_\theta(x,y) u\in C(\mathbb{L}'),~~~D_\theta(x,y) \in \mathrm{Der}(\mathbb{L}'), \label{2.10}\\
& D_\theta(x,y) [u,v,w]'=0,~~ \text{for all}~~ x,y\in\mathbb{L}~\text{and}~~u,v,w\in\mathbb{L}'.\label{2.11}
\end{align}

\begin{exam}
Let $(\mathbb{L}, [-,-,-])$  be a Lie triple system. If $\mathbb{L}^1$ satisfies $\mathbb{L}^1\in C(\mathbb{L})$, then the adjoint
representation  $\theta: \mathbb{L}\times \mathbb{L}\rightarrow \mathrm{End}(\mathbb{L})$ is an (coherent) action of $\mathbb{L}$ on itself.
\end{exam}

\begin{definition} \cite{Casas,Kasymov}
A 3-Leibniz algebra is a vector space $\mathfrak{g}$ endowed with a   multilinear map $[-, -, -]_\mathfrak{g}:\mathfrak{g}\times\mathfrak{g}\times\mathfrak{g}\rightarrow\mathfrak{g}$ satisfying
\begin{align}
%&[x,y,z]+[y,x,z]=0, \nonumber\\
 &[s, t, [u, v, w]_\mathfrak{g}]_\mathfrak{g}=[[s, t, u]_\mathfrak{g},v,w]_\mathfrak{g}+ [u,  [s, t, v]_\mathfrak{g},w]_\mathfrak{g}+ [u,v,[s, t, w]_\mathfrak{g}]_\mathfrak{g},\label{2.12}
\end{align}
for all $ s, t, u, v, w\in \mathfrak{g}$.
\end{definition}

\begin{prop} \label{prop:nonabelian hemisemidirect product}
Let  $ \theta$  be a coherent action of  a  Lie triple system  $(\mathbb{L},[-,-,-])$ on another
 Lie triple system  $(\mathbb{L}',[-,-,-]')$.
Then $(\mathbb{L}\oplus \mathbb{L}',[-,-,-]_{\ltimes})$ is a  3-Leibniz algebra, where
\begin{align}
[(a,u),(b,v),(c,w)]_{\ltimes}&=([a,b,c],D_\theta(a,b)w+[u,v,w]'),\label{2.13}
\end{align}
for any $(a,u),(b,v),(c,w)\in \mathbb{L}\oplus \mathbb{L}'$.
This   3-Leibniz algebra is called the nonabelian hemisemidirect product  3-Leibniz algebra and denoted by  $ \mathbb{L}\ltimes \mathbb{L}'$.
\end{prop}

\begin{proof}
For any $(x,u),(y,v),(z,w),(a,s),(b,t)\in \mathbb{L}\oplus  \mathbb{L}'$, by Eqs.  \eqref{2.3}, \eqref{2.7}, \eqref{2.10} and \eqref{2.11}, we have
\begin{align*}
&[[(a,s),(b,t),(x,u)]_\ltimes,(y,v),(z,w)]_\ltimes+ [(x,u),  [(a,s),(b,t), (y,v)]_\ltimes,(z,w)]_\ltimes+ \\
&[(x,u),(y,v),[(a,s),(b,t), (z,w)]_\ltimes]_\ltimes-[(a,s),(b,t), [(x,u), (y,v),(z,w)]_\ltimes]_\ltimes\\
&=[([a,b,x],D_\theta(a,b)u+[s,t,u]'),(y,v),(z,w)]_\ltimes+ [(x,u),  ([a,b,y],D_\theta(a,b)v+[s,t,v]'),(z,w)]_\ltimes+ \\
&~~~~[(x,u),(y,v),([a,b,z],D_\theta(a,b)w+[s,t,w]')]_\ltimes-[(a,s),(b,t), ([x,y,z],D_\theta(x,y)w+[u,v,w]')]_\ltimes\\
&=\big([[a,b,x],y,z], D_\theta([a,b,x], y)w+[D_\theta(a,b)u+[s,t,u]',v,w]'\big)+ \\
&~~~~~~\big([x,[a,b,y],z], D_\theta(x,[a,b,y])w+[u,D_\theta(a,b)v+[s,t,v]',w]'\big)+ \\
&~~~~~~\big([x,y,[a,b,z]], D_\theta(x,y)(D_\theta(a,b)w+[s,t,w]')+[u,v,D_\theta(a,b)w+[s,t,w]']'\big)-\\
&~~~~~~\big([a,b,[x,y,z]], D_\theta(a,b)(D_\theta(x,y)w+[u,v,w]')+[s,t,D_\theta(x,y)w+[u,v,w]']'\big)\\
&=\Big([[a,b,x],y,z]+[x,[a,b,y],z]+[x,y,[a,b,z]]-[a,b,[x,y,z]],\\
&~~~~~~ D_\theta([a,b,x], y)w+D_\theta(x,[a,b,y])w+D_\theta(x,y)D_\theta(a,b)w-D_\theta(a,b)D_\theta(x,y)w+\\
&~~~~~~[D_\theta(a,b)u,v,w]'+[[s,t,u]',v,w]'+[u,D_\theta(a,b)v,w]'+[u, [s,t,v]',w]'+[u,v,D_\theta(a,b)w ]'+\\
&~~~~~~[u,v, [s,t,w]']'+D_\theta(x,y) [s,t,w]'-D_\theta(a,b) [u,v,w]'- [s,t,D_\theta(x,y)w ]'-[s,t, [u,v,w]']' \Big)\\
&=0.
\end{align*}
Therefore,  $ \mathbb{L}\ltimes \mathbb{L}'$  is a  3-Leibniz algebra.
\end{proof}

\begin{definition}
Let $(\mathbb{L},[-,-,-])$ and  $(\mathbb{L}',[-,-,-]')$ be two Lie triple systems.
Let  $ \theta$  be a coherent action of   $\mathbb{L}$ on $\mathbb{L}'$.
   A linear map $T:\mathbb{L}'\rightarrow \mathbb{L}$ is called a nonabelian embedding tensor
   from $\mathbb{L}'$ to $\mathbb{L}$
with respect to     $\theta$   if $T$ satisfies
\begin{align}
 [Tu,Tv,Tw]=T(D_\theta(Tu,Tv)w+[u,v,w]'),\label{2.14}
\end{align}
for $u,v,w\in \mathbb{L}'$.
\end{definition}

\begin{remark}
If   $(\mathbb{L}', [-,-,-]')$   is an abelian  Lie triple system, then we can get that $T$ is an  embedding tensor on a Lie triple system  $(\mathbb{L},[-,-,-])$ with respect to the representation $(\mathbb{L}',\theta)$,  see \cite{Teng25} for more details. In addition, If $\theta=0$, then $T$ is a Lie triple system homomorphism from $(\mathbb{L}', [-,-,-]')$ to $(\mathbb{L}, [-,-,-])$.
\end{remark}

\begin{definition}
  Let $T$ and $\widetilde{T}$
 be two nonabelian embedding tensors
 from a Lie triple system  $(\mathbb{L}',[-,-,-]')$  to another Lie triple system  $(\mathbb{L},[-,-,-])$ with respect to  a coherent action $\theta$. A homomorphism from   $\widetilde{T}$
to $T$ consists of two Lie triple system  homomorphisms $f:\mathbb{L}\rightarrow\mathbb{L} $ and $f':\mathbb{L}'\rightarrow\mathbb{L}' $ such that
\begin{align}
  f(\widetilde{T}u) =&T  f'(u),\label{2.15}\\
 f'(\theta (x,y)u)=&\theta(f(x),f(y))f'(u),  \label{2.16}
\end{align}
for all $x,y\in\mathbb{L}$ and $u\in \mathbb{L}'.$
In particular, if $f$ and $f'$ are invertible, $(f,f')$ is said to be an isomorphism from  $\widetilde{T}$  to  $T$.
\end{definition}

 Let  $(f,f')$ is a homomorphism from  $\widetilde{T}$ to $T$, then we have
\begin{align}
 f'(D_\theta (x,y)u)=&D_\theta(f(x),f(y))f'(u),\label{2.17}
\end{align}

\begin{prop} \label{prop:newnonabelian}
Let $T:\mathbb{L}'\rightarrow \mathbb{L}$
 be a nonabelian embedding tensor
 from a Lie triple system  $(\mathbb{L}',[-,-,-]')$  to another Lie triple system  $(\mathbb{L},[-,-,-])$ with respect to  a coherent action $\theta$.
  Let $f:\mathbb{L}\rightarrow\mathbb{L} $ and $f':\mathbb{L}'\rightarrow\mathbb{L}' $ be two Lie triple system  isomorphisms  such that  Eqs. \eqref{2.15}
and \eqref{2.16} hold. Then $f^{-1}\circ T\circ f':\mathbb{L}'\rightarrow \mathbb{L}$ is a nonabelian embedding tensor from a Lie triple system  $(\mathbb{L}',[-,-,-]')$
to another Lie triple system  $(\mathbb{L},[-,-,-])$ with respect to  a coherent action $\theta$.
\end{prop}

\begin{proof}
For any $u,v,w\in \mathbb{L}'$,  we have
\begin{align*}
&[(f^{-1}\circ T\circ f')(u),(f^{-1}\circ T\circ f')(v),(f^{-1}\circ T\circ f')(w)]\\
&=f^{-1}[T  f'(u), T  f'(v), T  f'(w)]\\
&=f^{-1}\big(T(D_\theta(Tf'(u),Tf'(v))f'(w)+[f'(u),f'(v),f'(w)]')\big)\\
&=T\big(f'^{-1}(D_\theta(Tf'(u),Tf'(v))f'(w))+f'^{-1}[f'(u),f'(v),f'(w)]'\big)\\
&=T\big(D_\theta(f^{-1}(Tf'(u)),f^{-1}(Tf'(v))) w+[u, v,w]'\big),
\end{align*}
which implies that   $f^{-1}\circ T\circ f'$ is a  nonabelian embedding tensor.
\end{proof}

Using the   Eq. \eqref{2.13} of nonabelian hemisemidirect product, one can characterize $T$ by their graphs.

\begin{prop}
A linear map $T:\mathbb{L}'\rightarrow \mathbb{L}$ is   a nonabelian embedding tensor from a Lie triple system  $(\mathbb{L}',[-,-,-]')$  to another Lie triple system  $(\mathbb{L},[-,-,-])$ with respect to  a coherent action $\theta$ if and only if the graph
  $$Gr(T)=\{(Tu,u)~|~u\in \mathbb{L}'\}$$
  is a subalgebra of the nonabelian  hemisemidirect product  3-Leibniz algebra $ \mathbb{L}\ltimes \mathbb{L}'$.
\end{prop}

\begin{proof}
Let $T:\mathbb{L}'\rightarrow \mathbb{L}$ be a linear map.
For any $u,v,w\in \mathbb{L}'$,  we have
\begin{align*}
&[(Tu,u),(Tv,v),(Tw,w)]_\ltimes =\big([Tu,Tv,Tw],D_\theta(Tu,Tv)w+[u,v,w]'\big).
\end{align*}
Therefore,  the graph $Gr(T)=\{(Tu,u)~|~u\in \mathbb{L}'\}$ is a subalgebra of the nonabelian hemisemidirect product  3-Leibniz algebra $ \mathbb{L}\ltimes \mathbb{L}'$ if and only if $T$ satisfies Eq. \eqref{2.14}, which implies
that $T$ is a  nonabelian embedding tensor.
\end{proof}

Since $Gr(T) $ is isomorphic to $\mathbb{L}'$ as a vector space. we get the following conclusion immediately.

\begin{coro} \label{coro:NET}
Let $T:\mathbb{L}'\rightarrow \mathbb{L}$ be   a nonabelian embedding tensor from a Lie triple system  $(\mathbb{L}',[-,-,-]')$  to another Lie triple system  $(\mathbb{L},[-,-,-])$ with respect to  a coherent action $\theta$.  Define a linear map $[-,-,-]_T: \mathbb{L}'\times \mathbb{L}'\times \mathbb{L}'\rightarrow \mathbb{L}'$ by
\begin{align}
 [u,v,w]_T=D_\theta(Tu,Tv)w+[u,v,w]',\label{2.18}
\end{align}
for $u,v,w\in \mathbb{L}'$.  Then $(\mathbb{L}',[-,-,-]_T)$ is a  3-Leibniz algebra. Moreover, $T$ is a homomorphism from the 3-Leibniz algebra $(\mathbb{L}',[-,-,-]_T)$ to the Lie triple system $(\mathbb{L},[-,-,-])$.
The 3-Leibniz algebra $(\mathbb{L}',[-,-,-]_T)$  is called the descendent 3-Leibniz algebra of the nonabelian embedding tensor $T$.
\end{coro}

\begin{exam} \label{exam:3-dimensional}
Let $(\mathbb{L}_3,[-,-,-])$ be a 3-dimensional Lie triple system with a basis $\{\epsilon_1,\epsilon_2,\epsilon_3\}$
and the nonzero multiplication is defined by
$$[\epsilon_1,\epsilon_2,\epsilon_1]=-[\epsilon_2,\epsilon_1,\epsilon_1]=\epsilon_3.$$
The center of $(\mathbb{L}_3,[-,-,-])$ is the subspace generated by $\{\epsilon_3\}$. It is obvious that the adjoint
representation $\theta: \mathbb{L}_3\times \mathbb{L}_3\rightarrow \mathrm{End}(\mathbb{L}_3)$ is a coherent action of $\mathbb{L}_3$ on itself.
Consider a linear transformation $T$ given by
$T=\begin{bmatrix}
 a_{1} & b_{1} & c_{1}\\
 a_{2} & b_{2} & c_{2}\\
 a_{3} & b_{3} & c_{3}
 \end{bmatrix}$
with respect to the basis  $\{\epsilon_1,\epsilon_2,\epsilon_3\}$.
It straightforward to check that $T$  is a nonabelian embedding tensor  from $(\mathbb{L}_3,[-,-,-])$ to $(\mathbb{L}_3,[-,-,-])$ with respect to
the adjoint action $\theta$ if and only if
\begin{align*}
 [T\epsilon_i,T\epsilon_j,T\epsilon_k]=T([T\epsilon_i,T\epsilon_j,\epsilon_k]+[\epsilon_i,\epsilon_j,\epsilon_k]),~~~i,j,k=1,2,3.
\end{align*}
Therefore,  $T$  is a nonabelian embedding tensor if and
only if $c_1=c_2=0$ and $a_1^2b_2-a_1a_2b_1=c_3(a_1b_2-a_2b_1+1)$.
In particular,  $T=\begin{bmatrix}
 2 &  0 & 0\\
 0  & \frac{1}{2} & 0\\
 0 & 0 & 1
 \end{bmatrix}$ is a  nonabelian embedding tensor.
 By Corollary \ref{coro:NET}, we can get the 3-dimensional descendent
3-Leibniz algebra $(\mathbb{L}_3,[-,-,-]_T)$ given with respect to a basis $\{\epsilon_1,\epsilon_2,\epsilon_3\}$ by
$$[\epsilon_1,\epsilon_2,\epsilon_1]_T=2\epsilon_3,~~[\epsilon_2,\epsilon_1,\epsilon_1]_T=-2\epsilon_3.$$
\end{exam}

\begin{prop}\label{pro:2.9}
Let $T$ and $\widetilde{T}$
 be two nonabelian embedding tensors from a Lie triple system  $(\mathbb{L}',[-,-,-]')$  to another Lie triple system  $(\mathbb{L},[-,-,-])$ with respect to  a coherent action $\theta$,
  and
  $(f, f')$    a homomorphism from $\widetilde{T}$ to $T$. Then $f'$  is a
homomorphism of descendent 3-Leibniz algebras  from $(\mathbb{L}',[-,-,-]_{\widetilde{T}})$ to $(\mathbb{L}',[-,-,-]_{T})$.
\end{prop}

\begin{proof}  For all $u,v,w\in \mathbb{L}',$  by Eqs.   \eqref{2.15}--\eqref{2.18}, we have
\begin{align*}
f'([u,v,w]_{ \widetilde{T}})=&f'(D_\theta(\widetilde{T}u,\widetilde{T}v)w+[u,v,w]')\\
=&D_\theta(f(\widetilde{T}u),f(\widetilde{T}v))f'(w)+[f'(u),f'(v),f'(w)]'\\
=&D_\theta(Tf'(u),Tf'(v))f'(w)+[f'(u),f'(v),f'(w)]'\\
=&[f'(u),f'(v),f'(w)]_{T}.
\end{align*}
 The proof is finished.
 \end{proof}

\section{ Maurer-Cartan characterization  of    nonabelian  embedding tensors on Lie triple systems }\label{sec: Maurer-Cartan}
\def\theequation{\arabic{section}.\arabic{equation}}
\setcounter{equation} {0}

In this section,
%we first recall the concept of  an $L_{\infty}$-algebra and an Lie 3-algebra. Subsequently,
we construct a suitable $L_{\infty}$-algebra whose Maurer-Cartan elements are  nonabelian  embedding tensors   on
Lie triple systems.   This characterization
allows us to define a cohomology associated with a  nonabelian  embedding tensor.  In that follow section, we interpret this as the  cohomology of the corresponding
descendent 3-Leibniz algebras with coefficients in a suitable representation.
 %Following this, we  establish the cohomology of nonabelian  embedding tensors on Lie triple systems and realize it as the cohomology of the descendent 3-Leibniz algebra with coefficients in a suitable
%representation.

 \begin{definition} \cite{Stasheff}
An $L_{\infty}$-algebra is a $ \mathbb{Z}$-graded vector space $L=\oplus_{k\in  \mathbb{Z}} L^k$
equipped with a collection $(k\geq 1)$  of linear maps $l_k: \otimes^k L\rightarrow L$ of degree 1 with the property that, for any homogeneous elements
$x_1,x_2,\cdots, x_n\in L$,

(i) (graded symmetry) for every $\sigma\in S_n$,
$$l_n(x_{\sigma(1)},x_{\sigma(2)},\ldots, x_{\sigma(n)})=\varepsilon(\sigma)l_n(x_1,x_2,\ldots, x_n).$$

(ii) (generalized Jacobi Identity) for all $n\geq 1$,
$$\sum_{i=1}^n\sum_{\sigma\in S_n}\varepsilon(\sigma)l_{n-i+1}(l_n(x_{\sigma(1)},x_{\sigma(2)},\ldots, x_{\sigma(i)}),x_{i+1},\ldots, x_n)=0.$$
\end{definition}

%The notion of a Lie 3-algebra is a special case of $L_{\infty}$-algebra, in which only the 3-ary bracket is nonzero.
 %\begin{definition} \cite{Getzler,Tang}
%A Lie 3-algebra is a   $ \mathbb{Z}$-graded vector space $L=\oplus_{k\in  \mathbb{Z}} L^k$
%equipped with a trilinear bracket $[\![-,-,-]\!]_L: \otimes^3 L\rightarrow L$ of degree 1 satisfying the following.
%
%(i) (graded symmetry) For all homogeneous elements $x_1,x_2, x_3\in L$,
%$$[\![x_1,x_2, x_3]\!]_L=(-1)^{x_1x_2}[\![x_2,x_1, x_3]\!]_L=(-1)^{x_2x_3}[\![x_1,x_3, x_2]\!]_L.$$
%
%(ii) (generalized Jacobi Identity) For all homogeneous elements $x_1,x_2,\cdots, x_5\in L$,
%$$\sum_{\sigma\in S_5}\varepsilon(\sigma)[\![[\![x_{\sigma(1)},x_{\sigma(2)}, x_{\sigma(3)}]\!]_L,x_{4},  x_5]\!]_L=0.$$
%\end{definition}

 \begin{definition} \cite{Getzler}
%(i)
 A Maurer-Cartan element of an $L_{\infty}$-algebra $(L=\oplus_{k\in  \mathbb{Z}} L^k,\{l_i\}_{i=1}^\infty)$ is an element $\alpha\in L^0$
satisfying the Maurer-Cartan equation,
$$\sum_{n=1}^{\infty}\frac{1}{n!}l_n(\alpha,\ldots,\alpha)=0.$$

%(ii) A Maurer-Cartan element of a Lie 3-algebra $(L=\oplus_{k\in  \mathbb{Z}} L^k,[\![-,-,-]\!]_L)$ is an element $\alpha\in L^0$
%satisfying the Maurer-Cartan equation,
%$$\frac{1}{3!}[\![\alpha,\alpha,\alpha]\!]_L=0.$$
\end{definition}

\begin{prop}\cite{Getzler}   \label{prop:Lie 3-algebra}
Let $\alpha$ be a Maurer-Cartan element of an $L_{\infty}$-algebra $(L,\{l_i\}_{i=1}^\infty)$. For all $k\geq 1$ and $x_1,x_2,  x_3\in  L$, define
$l_k^\alpha: \otimes^k L\rightarrow L$  of degree 1 by
\begin{align*}
l_k^\alpha(x_1,x_2,\ldots, x_k)=&\sum_{n=0}^{\infty}\frac{1}{n!}l_{n+k}(\underbrace{\alpha,\ldots,\alpha}_n,x_1,x_2,\ldots, x_k).
%l_2^\alpha(x_1,x_2)=&[\![\alpha,x_1,x_2]\!]_L,\\
%l_3^\alpha(x_1,x_2,x_3)=&[\![x_1,x_2,x_3]\!]_L,\\
%l_k^\alpha=&0, k\geq 4.
\end{align*}
Then, $(L,\{l_k^\alpha\}_{k=1}^{+\infty})$ is an $L_\infty$-algebra, obtained from the
$L_{\infty}$-algebra  $(L,\{l_i\}_{i=1}^\infty)$ by twisting with the Maurer-Cartan element $\alpha$.  Moreover, $\alpha+\alpha'$
is a
Maurer-Cartan element of $(L,\{l_i\}_{i=1}^{+\infty})$  if and only if $\alpha'$
is a Maurer-Cartan element of the twisted $L_\infty$-algebra $(L,\{l_k^\alpha\}_{k=1}^{+\infty})$.
\end{prop}

In the sequel, we recall  Voronov's derived brackets theory \cite{Voronov}, which is a useful tool to
construct explicit  $L_\infty$-algebras.

 \begin{definition} \cite{Voronov}
A $V$-data consists of a quadruple $( \mathfrak{L}, \mathfrak{F},\mathcal{P},\Delta)$, where

$\bullet$  $(\mathfrak{L}, [-,-])$ is a graded Lie algebra;

$\bullet$  $\mathfrak{F}$ is an abelian graded Lie subalgebra of $(\mathfrak{L}, [-,-])$;

$\bullet$  $\mathcal{P}:\mathfrak{L}\rightarrow\mathfrak{L}$ is a projection, that is, $\mathcal{P}\circ\mathcal{P}=\mathcal{P}$, whose image is  $\mathfrak{F}$ and kernel is a graded Lie
subalgebra of $(\mathfrak{L}, [-,-])$;

$\bullet$  $\Delta$ is an element in $\mathrm{ker}(\mathcal{P})^1$ such that $[\Delta,\Delta]=0$.
 \end{definition}

\begin{prop}\cite{Voronov} \label{prop:V-data}
Let $(\mathfrak{L},\mathfrak{F},\mathcal{P},\Delta)$  be a $V$-data. Then $(\mathfrak{F},\{l_i\}_{i=1}^\infty)$ is an $L_\infty$-algebra, where
\begin{align}
 l_i(x_1,x_2,\ldots,x_i)=\mathcal{P}\underbrace{[\cdots[[}_i\Delta,x_1],x_2],\ldots,x_i]  ,\label{3.0}
 \end{align}
for homogeneous $x_1,x_2,\ldots,x_i\in \mathfrak{F}$.  We call $\{l_i\}_{i=1}^\infty$
the higher derived brackets of the $V$-data $(\mathfrak{L},\mathfrak{F},\mathcal{P},\Delta)$ .
\end{prop}

Let $L$ be a vector space. We consider the graded vector space $\mathfrak{C}_{\mathrm{3LA}}^\ast(L,L)=\oplus_{n\geq -1}\mathfrak{C}_{\mathrm{3LA}}^{n+1}(L,L)$, where
$\mathfrak{C}_{\mathrm{3LA}}^0(L,L)=\wedge^2 L$ and
 $\mathfrak{C}_{\mathrm{3LA}}^{n+1}(L,L)$
is the set of linear maps
$$P\in \mathrm{Hom}\big(\underbrace{(\wedge^2 L)\otimes\cdots\otimes(\wedge^2 L)}_n\otimes L,L\big).$$
 %satisfying
%$$P(\mathcal{X}_1,\ldots,\mathcal{X}_{n-1},x,x,y)=0,$$
%for $\mathcal{X}_i\in \otimes^2 g, 1\leq i\leq n-1$ and $  x,y\in g.$

\begin{prop}\cite{Voronov}
The graded vector space $\mathfrak{C}_{\mathrm{3LA}}^\ast(L,L)$ equipped with the graded
bracket
\begin{align*}
[P,Q]_{\mathrm{3LA}}=P\circ Q-(-1)^{pq}Q\circ P,~~~ \forall P\in \mathfrak{C}_{\mathrm{3LA}}^{p+1}(L,L),Q\in \mathfrak{C}_{\mathrm{3LA}}^{q+1}(L,L),
\end{align*}
 is a graded Lie algebra.
where $P\circ Q\in \mathfrak{C}_{\mathrm{3LA}}^{p+q+1}(L,L)$ is defined by
\begin{align*}
&(P\circ Q)(\mathcal{X}_1,\ldots,\mathcal{X}_{p+q},x)\\
=&\sum_{k=1}^p(-1)^{(k-1)q}\sum_{\sigma\in S(k-1,q)}(-1)^\sigma P\Big(\mathcal{X}_{\sigma(1)},\ldots,\mathcal{X}_{\sigma{(k-1)}},Q\big(\mathcal{X}_{\sigma(k)},\ldots,\mathcal{X}_{\sigma{(k+q-1)}},x_{(k+q)}\big)\otimes y_{k+q},\\
&\mathcal{X}_{k+q+1},\ldots,\mathcal{X}_{p+q},x\Big)+\sum_{k=1}^p(-1)^{(k-1)q}\sum_{\sigma\in S(k-1,q)}(-1)^\sigma P\Big(\mathcal{X}_{\sigma(1)},\ldots,\mathcal{X}_{\sigma{(k-1)}}, \\
&x_{k+q}\otimes Q\big(\mathcal{X}_{\sigma(k)},\ldots,\mathcal{X}_{\sigma{(k+q-1)}},y_{k+q}\big),\mathcal{X}_{k+q+1},\ldots,\mathcal{X}_{p+q},x\Big)\\
&+\sum_{\sigma\in S(p,q)}(-1)^{pq}(-1)^\sigma P\Big(\mathcal{X}_{\sigma(1)},\ldots,\mathcal{X}_{\sigma{(p)}},Q\big(\mathcal{X}_{\sigma(p+1)},\ldots,\mathcal{X}_{\sigma{(p+q)}},x\big)\Big),
\end{align*}
for all $\mathcal{X}_i=x_i\wedge y_i\in \wedge^2 L, i=1,2,\ldots,p+q$ and $x\in L$.
\end{prop}

\begin{lemma}  \label{lemma:NA structure}
For $\pi\in \mathfrak{C}_{\mathrm{3LA}}^2(L,L)$, we have
\begin{align*}
&[\pi,\pi]_{\mathrm{3LA}}(\mathcal{X}_1, \mathcal{X}_{2},x)=2(\pi\circ\pi)(\mathcal{X}_1, \mathcal{X}_{2},x)\\
=&2\Big(\pi(\pi(x_1,y_1,x_2),y_2,x)+\pi(x_2,\pi(x_1,y_1,y_2), x)+\pi(x_2,y_2,\pi(x_1,y_1,x))-\pi(x_1,y_1,\pi(x_2,y_2,x))\Big).
\end{align*}
Thus, $\pi$ is a 3-Leibniz algebras structure if and only if $[\pi,\pi]_{\mathrm{3LA}}=0$, i.e. $\pi$  is a Maurer-Cartan
element of the graded Lie algebra $(\mathfrak{C}_{\mathrm{3LA}}^\ast(L,L),[-,-]_{\mathrm{3LA}})$.
\end{lemma}

Let  $ \theta$  be a coherent action of  a  Lie triple system  $(\mathbb{L},[-,-,-])$ on another
 Lie triple system  $(\mathbb{L}',[-,-,-]')$.
  Usually we use $\varpi$ to indicate the Lie triple system bracket  $[-,-,-]$ and $\varpi'$ to indicate the Lie triple system bracket  $[-,-,-]'$.
  Then the  nonabelian hemisemidirect product 3-Leibniz algebra given by Proposition \ref{prop:nonabelian hemisemidirect product} corresponds to
%  Then, $\varpi\boxtimes D_\theta\boxtimes \varpi'$ corresponds to the nonabelian hemisemidirect product 3-Leibniz algebra  structure on $\mathbb{L} \oplus  \mathbb{L}'$ given by
$$\varpi\boxtimes D_\theta \boxtimes \varpi'\big((a,u),(b,v),(c,w)\big)=([a,b,c],D_\theta(a,b)w+[u,v,w]'),$$
for  $(a,u),(b,v),(c,w)\in \mathbb{L}\oplus \mathbb{L}'$.  Therefore, by Lemma \ref{lemma:NA structure}, we have
$$[\varpi\boxtimes D_\theta\boxtimes \varpi',\varpi\boxtimes D_\theta\boxtimes \varpi']_{\mathrm{3LA}}=0.$$

\begin{lemma} \label{lemma:ETLTS}
Let $T:\mathbb{L}'\rightarrow \mathbb{L}$ be   a nonabelian embedding tensor from a Lie triple system  $(\mathbb{L}',[-,-,-]')$  to another Lie triple system  $(\mathbb{L},[-,-,-])$ with respect to  a coherent action $\theta$.  Then, we have
\begin{align*}
&[\varpi\boxtimes D_\theta \boxtimes \varpi', T]_{\mathrm{3LA}}\big((a,u),(b,v),(c,w)\big)\\
=&\big([Tu,b,c]+[a,b,Tw]+[a,Tv,c]-T(D_\theta(a,b)w+[u,v,w]'),D_\theta(a,Tv)w+D_\theta(Tu,b)w  \big);\\
&[[\varpi\boxtimes D_\theta \boxtimes \varpi', T]_{\mathrm{3LA}}, T]_{\mathrm{3LA}}\big((a,u),(b,v),(c,w)\big)\\
=&2\big([Tu,Tv,c]+[Tu,b,Tw]+[a,Tv,Tw]-T(D_\theta(Tu,b)w+D_\theta(a,Tv)w),D_\theta(Tu,Tv)w  \big),
\end{align*}
for all $(a,u),(b,v),(c,w)\in \mathbb{L}\oplus  \mathbb{L}'$.
\end{lemma}

\begin{proof}
For any  $(a,u),(b,v),(c,w)\in \mathbb{L}\oplus  \mathbb{L}'$, we have
\begin{align*}
&[\varpi\boxtimes D_\theta\boxtimes \varpi', T]_{\mathrm{3LA}}\big((a,u),(b,v),(c,w)\big)\\
&=((\varpi\boxtimes D_\theta\boxtimes \varpi')\circ T- T\circ (\varpi\boxtimes D_\theta\boxtimes \varpi'))\big((a,u),(b,v),(c,w)\big)\\
&=\varpi\boxtimes D_\theta\boxtimes \varpi'\big(T(a,u),(b,v),(c,w)\big)+\varpi\boxtimes D_\theta\boxtimes \varpi'\big((a,u),T(b,v),(c,w)\big)+\\
&~~~~~\varpi\boxtimes D_\theta\boxtimes \varpi'\big((a,u),(b,v),T(c,w)\big)-T([a,b,c],D_\theta(a,b)w+[u,v,w]')\\
&=\big([Tu,b,c]+[a,b,Tw]+[a,Tv,c]-T(D_\theta(a,b)w+[u,v,w]'),D_\theta(a,Tv)w+D_\theta(Tu,b)w  \big),\\
&[[\varpi\boxtimes D_\theta \boxtimes \varpi', T]_{\mathrm{3LA}}, T]_{\mathrm{3LA}}\big((a,u),(b,v),(c,w)\big)\\
&=[\varpi\boxtimes D_\theta\boxtimes \varpi', T]_{\mathrm{3LA}}\big(T(a,u),(b,v),(c,w)\big)+[\varpi\boxtimes D_\theta\boxtimes \varpi', T]_{\mathrm{3LA}}\big((a,u),T(b,v),(c,w)\big)+\\
&~~~~~[\varpi\boxtimes D_\theta\boxtimes \varpi', T]_{\mathrm{3LA}}\big((a,u),(b,v),T(c,w)\big)-T[\varpi\boxtimes D_\theta\boxtimes \varpi', T]_{\mathrm{3LA}}\big((a,u),(b,v),(c,w)\big)\\
%&=\big([[Tu,b,Tw]+[Tu,Tv,c],D_\theta(Tu,Tv)w  \big)+\big([Tu,Tv,c]+[a,Tv,Tw], D_\theta(Tu,Tv)w  \big)+\\
%&\big([Tu,b,Tw]+[a,Tv,Tw],0\big)-\big(TD_\theta(a,Tv)w+TD_\theta(Tu,b)w,0\big)\\
&=2\big([Tu,Tv,c]+[Tu,b,Tw]+[a,Tv,Tw]-T(D_\theta(Tu,b)w+D_\theta(a,Tv)w),D(Tu,Tv)w  \big).
\end{align*}
The proof is finished.
\end{proof}

\begin{prop}
Let  $ \theta$  be a coherent action of  a  Lie triple system  $(\mathbb{L},[-,-,-])$ on another
 Lie triple system  $(\mathbb{L}',[-,-,-]')$.
 Then we have a
$V$-data $( \mathfrak{L}, \mathfrak{F},\mathcal{P},\Delta)$ as follows:

$\bullet$  the graded Lie algebra $(\mathfrak{L},[-,-])$ is given by $(\mathfrak{C}_{\mathrm{3LA}}^\ast(\mathbb{L}\oplus  \mathbb{L}',\mathbb{L}\oplus  \mathbb{L}'),[-,-]_{\mathrm{3LA}})\mathrm{;}$

$\bullet$   the abelian graded Lie subalgebra $ \mathfrak{F}$ is given by
$$\mathfrak{F}=\mathfrak{C}^\ast( \mathbb{L}', \mathbb{L})=\oplus_{n=0}^{+\infty}\mathfrak{C}^{n+1}( \mathbb{L}',\mathbb{L})=\oplus_{n=0}^{+\infty}\mathrm{Hom}\big(\underbrace{(\wedge^2  \mathbb{L}')\otimes\cdots\otimes(\wedge^2  \mathbb{L}')}_{n\geq 0}\otimes  \mathbb{L}',\mathbb{L}\big);$$

$\bullet$  $\mathcal{P}: \mathfrak{L}\rightarrow \mathfrak{L}$ is the projection onto the subspace $ \mathfrak{F}$;

$\bullet$  $\Delta=\varpi\boxtimes D_\theta\boxtimes  \varpi'$.\\
Consequently, we obtain an  $L_\infty$-algebra $(\mathfrak{C}^\ast( \mathbb{L}',\mathbb{L}),l_1,l_3)$, where
\begin{align}
l_1(P)&=[\varpi\boxtimes D_\theta\boxtimes \varpi',P]_{\mathrm{3LA}}\in \mathfrak{C}^{m+1}(\mathbb{L}',\mathbb{L}),\\
l_3(P,Q,R)&=[[[\varpi\boxtimes D_\theta\boxtimes \varpi',P]_{\mathrm{3LA}},Q]_{\mathrm{3LA}},R]_{\mathrm{3LA}}\in \mathfrak{C}^{m+n+k+2}( \mathbb{L}',\mathbb{L}),\label{3.1}
\end{align}
for $P\in  \mathfrak{C}^{m+1}( \mathbb{L}',\mathbb{L})$, $Q\in  \mathfrak{C}^{n+1}( \mathbb{L}',\mathbb{L})$ and $R\in  \mathfrak{C}^{k+1}( \mathbb{L}',\mathbb{L})$.
\end{prop}

\begin{proof}
 It is obvious that
$\Delta=\varpi\boxtimes D_\theta\boxtimes  \varpi'\in \mathrm{ker}(\mathcal{P})^2$.
Therefore, $( \mathfrak{L}, \mathfrak{F},\mathcal{P},\Delta)$  is a $V$-data.
By Proposition \ref{prop:V-data}, $( \mathfrak{F},\{l_i\}_{i=1}^\infty)$ is an $L_\infty$-algebra, where $l_i$ is given by Eq. \eqref{3.0}. Moreover,  for any  $P\in  \mathfrak{C}^{m+1}(\mathbb{L}',\mathbb{L})$ and $Q\in  \mathfrak{C}^{n+1}(\mathbb{L}',\mathbb{L})$, by Lemma \ref{lemma:ETLTS}, we have
%$$[\varpi\boxtimes D_\theta\boxtimes  \varpi',P]_{\mathrm{3LA}}\in \mathrm{ker}(\mathcal{P}),$$
$$[[\varpi\boxtimes D_\theta\boxtimes  \varpi',P]_{\mathrm{3LA}},Q]_{\mathrm{3LA}}\in \mathrm{ker}(\mathcal{P}),$$
which implies that $l_2=0$. Similarly, we have $l_k=0$, when $k\geq 4$. Therefore,
 $(\mathfrak{C}^\ast( \mathbb{L}',\mathbb{L}),l_1,l_3)$  is an  $L_\infty$-algebra.
\end{proof}

\begin{theorem} \label{theorem:MC}
Let  $ \theta$  be a coherent action of  a  Lie triple system  $(\mathbb{L},[-,-,-])$ on another
 Lie triple system  $(\mathbb{L}',[-,-,-]')$. Then Maurer-Cartan
elements of the  $L_\infty$-algebra $(\mathfrak{C}^\ast(\mathbb{L}',\mathbb{L}),l_1,l_3)$ are precisely nonabelian embedding tensors from a Lie triple system  $(\mathbb{L}',[-,$ $-,-]')$  to another Lie triple system  $(\mathbb{L},[-,-,-])$ with respect to  a coherent action $\theta$.
\end{theorem}
\begin{proof}
For any $T\in \mathfrak{C}^1(\mathbb{L}',\mathbb{L})$ and $u,v,w\in \mathbb{L}'$. By Lemma \ref{lemma:ETLTS}, we can deduce that
\begin{align*}
&l_1(T)(u,v,w)=[\varpi\boxtimes D\boxtimes \varpi',T]_{\mathrm{3LA}}(u,v,w)=-T[u,v,w]',\\
&l_3(T,T,T)(u,v,w)\\
&=[[[\varpi\boxtimes D_\theta\boxtimes  \varpi',T]_{\mathrm{3LA}},T]_{\mathrm{3LA}},T]_{\mathrm{3LA}}(u,v,w)\\
&=[[\varpi\boxtimes D_\theta\boxtimes  \varpi',T]_{\mathrm{3LA}},T]_{\mathrm{3LA}}(Tu,v,w)+[[\varpi\boxtimes D_\theta\boxtimes  \varpi',T]_{\mathrm{3LA}},T]_{\mathrm{3LA}}(u,Tv,w)+\\
&~~~~~[[\varpi\boxtimes D_\theta\boxtimes  \varpi',T]_{\mathrm{3LA}},T]_{\mathrm{3LA}}(u,v,Tw)-T[[\varpi\boxtimes D_\theta\boxtimes  \varpi',T]_{\mathrm{3LA}},T]_{\mathrm{3LA}}(u,v,w)\\
%&=[\varpi\boxtimes D_\theta\boxtimes  \varpi',T]_{\mathrm{3LA}}(Tu,Tv,w)+[\varpi\boxtimes D_\theta\boxtimes  \varpi',T]_{\mathrm{3LA}}(Tu,v,Tw)-T[\varpi\boxtimes D_\theta\boxtimes  \varpi',T]_{\mathrm{3LA}}(Tu,v,w)+\\
%&~~~~~[\varpi\boxtimes D_\theta\boxtimes  \varpi',T]_{\mathrm{3LA}}(Tu,Tv,w)+[\varpi\boxtimes D_\theta\boxtimes  \varpi',T]_{\mathrm{3LA}}(u,Tv,Tw)-T[\varpi\boxtimes D_\theta\boxtimes  \varpi',T]_{\mathrm{3LA}}(u,Tv,w)+\\
%&~~~~~[\varpi\boxtimes D_\theta\boxtimes  \varpi',T]_{\mathrm{3LA}}(Tu,v,Tw)+[\varpi\boxtimes D_\theta\boxtimes  \varpi',T]_{\mathrm{3LA}}(u,Tv,Tw)-T[\varpi\boxtimes D_\theta\boxtimes  \varpi',T]_{\mathrm{3LA}}(u,v,Tw)-\\
%&~~~~~T[\varpi\boxtimes D_\theta\boxtimes  \varpi',T]_{\mathrm{3LA}}(Tu,v,w)-T[\varpi\boxtimes D_\theta\boxtimes  \varpi',T]_{\mathrm{3LA}}(u,Tv,w)-T[\varpi\boxtimes D_\theta\boxtimes  \varpi',T]_{\mathrm{3LA}}(u,v,Tw)\\
%&=[Tu,Tv,Tw]+[Tu,Tv,Tw] -T D_\theta (Tu, Tv)w+[Tu,Tv,Tw]+[Tu,Tv,Tw]-T D_\theta (Tu, Tv)w+\\
%&~~~~~[Tu,Tv,Tw]+[Tu,Tv,Tw]-T[\varpi\boxtimes D_\theta\boxtimes  \varpi',T]_{\mathrm{3LA}}(u,v,Tw)-\\
%&~~~~~T D_\theta (Tu, Tv)w-T D_\theta (Tu, Tv)w-T[\varpi\boxtimes D_\theta\boxtimes  \varpi',T]_{\mathrm{3LA}}(u,v,Tw)\\
%&=2[Tu,Tv,Tw]-TD_\theta(Tu,Tv)w+2[Tu,Tv,Tw]-TD_\theta(Tu,Tv)w+\\
%&~~~~~2[Tu,Tv,Tw]-2TD_\theta(Tu,Tv)w\\
&=6([Tu,Tv,Tw]-TD_\theta(Tu,Tv)w),
\end{align*}
Let $T$ be a Maurer-Cartan element of the $L_\infty$-algebra $(\mathfrak{C}^\ast(\mathbb{L}',\mathbb{L}),l_1,l_3)$.  We have
\begin{align*}
&(l_1(T) +\frac{1}{3!}l_3(T,T,T))(u,v,w)\\
 &=[\varpi\boxtimes D\boxtimes \varpi',T]_{\mathrm{3LA}}(u,v,w)+\frac{1}{6}[[[\varpi\boxtimes D_\theta\boxtimes  \varpi',T]_{\mathrm{3LA}},T]_{\mathrm{3LA}},T]_{\mathrm{3LA}}(u,v,w)\\
 &=-T[u,v,w]'+[Tu,Tv,Tw]-TD_\theta(Tu,Tv)w\\
 &=0,
\end{align*}
which implies that $T$ is a   nonabelian embedding tensor.
\end{proof}

\begin{prop}
 Let $T:\mathbb{L}'\rightarrow \mathbb{L}$ be   a nonabelian embedding tensor from a Lie triple system  $(\mathbb{L}',[-,-,-]')$  to another Lie triple system  $(\mathbb{L},[-,-,-])$ with respect to  a coherent action $\theta$.
 Then   $(\mathfrak{C}^\ast(\mathbb{L}',\mathbb{L}),l_1^T,l_2^T,l_3^T)$ is a twisted $L_\infty$-algebra, where $l_1^T,l_2^T,l_3^T$ are given by
\begin{align*}
l_1^T(P)=&l_1(P)+\frac{1}{2}l_3(T,T,P),\\
l_2^T(P,Q)=&l_3(T,P,Q),\\
l_3^T(P,Q,R)=&l_3(P,Q,R),\\
l_k^T=&0, k\geq 4,
\end{align*}
for all $P\in  \mathfrak{C}^{p}(\mathbb{L}',\mathbb{L})$, $Q\in  \mathfrak{C}^{q}(\mathbb{L}',\mathbb{L})$ and $R\in  \mathfrak{C}^{r}(\mathbb{L}',\mathbb{L})$.
\end{prop}

\begin{proof}
 Since the nonabelian embedding tensor $T$ is a Maurer-Cartan element of the   $L_\infty$-algebra $(\mathfrak{C}^\ast(\mathbb{L}',\mathbb{L}),l_1,l_3)$,
 then by Proposition \ref{prop:Lie 3-algebra}, we have the conclusions.
\end{proof}

Actually, the above twisted   $L_\infty$-algebra controls deformations of nonabelian embedding tensor $T$.
\begin{theorem} \label{theorem:MCD}
Let $T:\mathbb{L}'\rightarrow \mathbb{L}$ be   a nonabelian embedding tensor from a Lie triple system  $(\mathbb{L}',[-,-,-]')$  to another Lie triple system  $(\mathbb{L},[-,-,-])$ with respect to  a coherent action $\theta$.
Then   for a   linear map $\widetilde{T}:\mathbb{L}'\rightarrow \mathbb{L},$
$T+\widetilde{T}$ is a  nonabelian embedding tensor
if and only if $\widetilde{T}$ is a Maurer-Cartan element of the twisted $L_\infty$-algebra $(\mathfrak{C}^\ast(\mathbb{L}',\mathbb{L}),l_1^T,l_2^T,l_3^T)$,
that is,
$\widetilde{T}$
satisfies the Maurer-Cartan equation:
$$l_1^T(\widetilde{T})+\frac{1}{2!}l_2^T(\widetilde{T},\widetilde{T})+\frac{1}{3!}l_3^T(\widetilde{T},\widetilde{T},\widetilde{T})=0.$$
\end{theorem}

\begin{proof}
By Theorem \ref{theorem:MC}, $T+\widetilde{T}$ is a nonabelian embedding tensor if and only if
$$l_1(T+\widetilde{T})+\frac{1}{3!}l_3(T+\widetilde{T},T+\widetilde{T},T+\widetilde{T})=0.$$
Applying $l_1(T)+\frac{1}{3!}l_3(T,T,T)=0$, the above equation is equivalent to
$$l_1(\widetilde{T})+\frac{1}{2}l_3(T,T,\widetilde{T})+\frac{1}{2}l_3(T,\widetilde{T},\widetilde{T})+\frac{1}{6}l_3(\widetilde{T},\widetilde{T},\widetilde{T})=0,$$
that is, $l_1^T(\widetilde{T})+\frac{1}{2!}l_2^T(\widetilde{T},\widetilde{T})+\frac{1}{3!}l_3^T(\widetilde{T},\widetilde{T},\widetilde{T})=0,$
 which implies that $\widetilde{T}$ is a Maurer-Cartan element of the twisted $L_\infty$-algebra  $(\mathfrak{C}^\ast(\mathbb{L}',\mathbb{L}),l_1^T,l_2^T,l_3^T)$.
\end{proof}

The above characterization of a nonabelian embedding tensor $T$ allows us to define a cohomology associated with $T$.
 %Therefore,   we consider the cochain complex: $(\mathfrak{C}_T^{\bullet}(\mathbb{L}',\mathbb{L}),\mathrm{d}_T^\bullet)$.
More precisely,
% \begin{align*}
 %\mathfrak{C}_E^{n}(\mathbb{L}',\mathbb{L}):=
% \begin{cases}
%\mathfrak{C}^{n}(\mathbb{L}',\mathbb{L}) ,&n\geq 1,\\
% \wedge^2 \mathbb{L},&n=0
%\end{cases}
% \end{align*}
define the differential operator $\mathrm{d}_T^n:\mathfrak{C}^{n}(\mathbb{L}',\mathbb{L})\rightarrow \mathfrak{C}^{n+1}(\mathbb{L}',\mathbb{L})$ by
\begin{align}
 \mathrm{d}_T^n(f)=l_1^T(f)=&l_1(f)+\frac{1}{2}l_3(T,T,f),~~ \forall ~f\in \mathfrak{C}_T^{n}(\mathbb{L}',\mathbb{L}), ~\text{for} ~~n\geq 1. \label{3.2}
\end{align}
%if $n=0,$  define $\mathrm{d}_T^0:\mathfrak{C}_T^{0}(\mathbb{L}',\mathbb{L})\rightarrow \mathfrak{C}_T^{1}(\mathbb{L}',\mathbb{L})$ by
%\begin{align}
% \mathrm{d}_T^0(a,b)u=T(D_\theta(a,b)u)-[a,b,Tu],~ \forall (a,b)\in \mathfrak{C}_T^{0}(\mathbb{L}',\mathbb{L}),~ u\in \mathbb{L}'. \label{3.3}
%\end{align}
 The corresponding $n$-th cohomology group of the  nonabelian embedding tensor $T$  is denoted by
\begin{align*}
 \mathcal{H}_T^n(\mathbb{L}',\mathbb{L})=\frac{\mathcal{Z}_T^n(\mathbb{L}',\mathbb{L})}{\mathcal{B}_T^n(\mathbb{L}',\mathbb{L})}=\frac{\{f\in \mathfrak{C}^{n}(\mathbb{L}',\mathbb{L})~|~ \mathrm{d}_T^n(f)=0\}}{\{\mathrm{d}_T^n(f)~|~f\in \mathfrak{C}^{n-1}(\mathbb{L}',\mathbb{L})\}}.
\end{align*}

\section{ Cohomology  of    nonabelian  embedding tensors between Lie triple systems }\label{sec: Cohomology}
\def\theequation{\arabic{section}.\arabic{equation}}
\setcounter{equation} {0}

In this section, we construct a representation of the descendent 3-Leibniz algebra $(\mathbb{L}',[-,-,-]_T)$
on the vector space $\mathbb{L}$ and define the cohomologies of a nonabelian  embedding tensor  between Lie triple systems.
For this purpose,    first we recall some basic results about  representations and cohomologies of
 3-Leibniz algebra $(\mathfrak{g},[-,-,-]_\mathfrak{g})$.

\begin{definition} \cite{Casas} A representation of the 3-Leibniz algebra $(\mathfrak{g},[-,-,-]_{\mathfrak{g}})$ is a
vector space $\mathcal{V}$ equipped with 3 actions
\begin{align*}
\rho^l:\mathfrak{g}\otimes \mathfrak{g}\otimes \mathcal{V}\rightarrow \mathcal{V},\\
\rho^m:\mathfrak{g}\otimes \mathcal{V}\otimes \mathfrak{g}\rightarrow \mathcal{V},\\
\rho^r:\mathcal{V}\otimes \mathfrak{g}\otimes \mathfrak{g}\rightarrow \mathcal{V},
\end{align*}
such that for any  $x,y,z,a,b\in \mathfrak{g}$ and $u\in \mathcal{V}$,
\begin{align}
&\rho^l(a,b,\rho^l(x,y,u))=\rho^l([a,b,x]_{\mathfrak{g}}, y, u)+ \rho^l(x,[a,b,y]_{\mathfrak{g}},u)+\rho^l(x,y,\rho^l(a,b,u)),\label{4.1}\\
 &\rho^l(a,b,\rho^m(x,u,z))= \rho^m([a,b,x]_{\mathfrak{g}}, u, z)+\rho^m(x,\rho^l(a,b,u),z)+\rho^m(x,u,[a,b,z]_{\mathfrak{g}}),\label{4.2}\\
  &\rho^l(a,b,\rho^r(u,y,z))=\rho^r(\rho^l(a,b,u), y,z)+ \rho^r(u,[a,b,y]_{\mathfrak{g}},z)+\rho^r(u,y,[a,b,z]_{\mathfrak{g}}),\label{4.3}\\
   &\rho^m(a,u,[x,y,z]_{\mathfrak{g}})=\rho^r(\rho^m(a,u,x), y,z)+ \rho^m(x,\rho^m(a,u,y),z)+\rho^l(x,y,\rho^m(a,u,z)),\label{4.4}\\
    &\rho^r(u,b,[x,y,z]_{\mathfrak{g}})=\rho^r(\rho^r(u,b,x), y,z)+ \rho^m(x,\rho^r(u,b,y),z)+\rho^l(x,y,\rho^r(u,b,z)).\label{4.5}
\end{align}
\end{definition}

An $n$-cochain on a 3-Leibniz algebra $(\mathfrak{g},[-,-,-]_{\mathfrak{g}})$  with coefficients in a representation  $(\mathcal{V};\rho^l,\rho^m,\rho^r)$
is a linear map
\begin{align*}
f:\underbrace{\wedge^2\mathfrak{g}\otimes\cdots \otimes \wedge^2 \mathfrak{g}}_{n-1}\otimes \mathfrak{g}\rightarrow \mathcal{V},n\geq 1.
\end{align*}
The space generated by $n$-cochains is denoted as $\mathcal{C}^n_{\mathrm{3Leib}}(\mathfrak{g},\mathcal{V})$.  The coboundary map $\delta^n: \mathcal{C}^n_{\mathrm{3Leib}}(\mathfrak{g},\mathcal{V})\rightarrow \mathcal{C}^{n+1}_{\mathrm{3Leib}}(\mathfrak{g},\mathcal{V})$, $f\mapsto (\delta^n f)$, for $\mathfrak{X}_i=x_i\wedge y_i\in \wedge^2 \mathfrak{g}(1\leq i\leq n)$ and $z\in \mathfrak{g}$, as
%\begin{small}
\begin{align*}
&(\delta^n f)(\mathfrak{X}_1,\mathfrak{X}_2, \ldots, \mathfrak{X}_n,z)\\
=&\sum_{1\leq j<k\leq n}(-1)^jf(\mathfrak{X}_1,\ldots,\widehat{\mathfrak{X}_j},\ldots,\mathfrak{X}_{k-1},x_k \wedge[x_j,y_j,y_k]_{\mathfrak{g}}+[x_j,y_j,x_k]_{\mathfrak{g}}\wedge y_k,\ldots,\mathfrak{X}_n,z)\\
&+\sum_{j=1}^n(-1)^jf(\mathfrak{X}_1,\cdots,\widehat{\mathfrak{X}_j},\ldots,\mathfrak{X}_{n},[x_j,y_j,z]_{\mathfrak{g}})+\sum_{j=1}^n(-1)^{j+1}\rho^l(\mathfrak{X}_j,f(\mathfrak{X}_1,\ldots,\widehat{\mathfrak{X}_j},\ldots,\mathfrak{X}_{n},z))\\
&+(-1)^{n+1}(\rho^m(x_n, f(\mathfrak{X}_1,\cdots,\mathfrak{X}_{n-1},y_n),z)+\rho^r(f(\mathfrak{X}_1,\cdots,\mathfrak{X}_{n-1},x_n), y_n,z)).
\end{align*}
%\end{small}
It was proved in \cite{Casas,Takhtajan} that $\delta^{n+1}\circ\delta^n=0$. Thus,
$(\oplus_{n=1}^{+\infty}\mathcal{C}^n_{\mathrm{3Leib}}(\mathfrak{g},\mathcal{V}),\delta^\bullet)$ is a cochain complex.  The corresponding  $n$-th cohomology group by $\mathcal{H}^n_{\mathrm{3Leib}}(\mathfrak{g},\mathcal{V})=\mathcal{Z}^n_{\mathrm{3Leib}}(\mathfrak{g},\mathcal{V})/\mathcal{B}^n_{\mathrm{3Leib}}(\mathfrak{g},\mathcal{V})$.

\begin{lemma} \label{lemma:representation}
 Let $T:\mathbb{L}'\rightarrow \mathbb{L}$ be   a nonabelian embedding tensor from a Lie triple system  $(\mathbb{L}',[-,-,-]')$  to another Lie triple system  $(\mathbb{L},[-,-,-])$ with respect to  a coherent action $\theta$. Define  actions
\begin{align*}
\rho^l_T:\mathbb{L}'\otimes \mathbb{L}'\otimes \mathbb{L}\rightarrow \mathbb{L},\rho^m_T:\mathbb{L}'\otimes \mathbb{L}\otimes \mathbb{L}'\rightarrow \mathbb{L},\rho^r_T:\mathbb{L}\otimes \mathbb{L}'\otimes \mathbb{L}'\rightarrow \mathbb{L},
\end{align*}
by
%\begin{small}
\begin{align}
\rho^l_T(u,v,x)&=[Tu,Tv,x],\label{4.6}\\
\rho^m_T(u,x,v)&=[Tu,x,Tv]-TD_\theta(Tu,x)v,\label{4.7}\\
\rho^r_T(x,u,v)&=[x,Tu,Tv]-TD_\theta(x,Tu)v,\label{4.8}
\end{align}
%\end{small}
for  $u,v\in \mathbb{L}'$ and $x\in \mathbb{L}.$ Then $(\mathbb{L};\rho^l_T,\rho^m_T,\rho^r_T)$  is a representation of the 3-Leibniz algebra $(\mathbb{L}',[-,-,-]_T)$  given  in
Corollary \ref{coro:NET}.
\end{lemma}

\begin{proof}
For any $u,v,s,t\in \mathbb{L}'$ and $x\in \mathbb{L}$, by Eqs. \eqref{2.3}, \eqref{2.14} and \eqref{2.18}, we have
\begin{align*}
&\rho^l_T(u,v,\rho^l_T(s,t,x))-\rho^l_T([u,v,s]_{T}, t, x)- \rho^l_T(s,[u,v,t]_{T},x)-\rho^l_T(s,t,\rho^l_T(u,v,x))\\
=&[Tu,Tv,[Ts,Tt,x]]-[T[u,v,s]_{T}, Tt, x]- [Ts,T[u,v,t]_{T},x]-[Ts,Tt,[Tu,Tv,x]]\\
=&[Tu,Tv,[Ts,Tt,x]]-[[Tu,Tv,Ts], Tt, x]- [Ts,[Tu,Tv,Tt],x]-[Ts,Tt,[Tu,Tv,x]]\\
=&0,
\end{align*}
which indicates that Eq. \eqref{4.1} follows.

Further, by Eqs. \eqref{2.3}, \eqref{2.7},  \eqref{2.14} and \eqref{2.17}, we have
\begin{align*}
 &\rho^l_T(u,v,\rho^m_T(s,x,t))- \rho^m_T([u,v,s]_{T}, x, t)-\rho^m_T(s,\rho^l_T(u,v,x),t)-\rho^m_T(s,x,[u,v,t]_{T})\\
 =&[Tu,Tv,[Ts,x,Tt]]-[Tu,Tv,TD_\theta(Ts,x)t]- [T[u,v,s]_{T}, x, Tt]+TD_\theta(T[u,v,s]_{T}, x)t\\
 &-[Ts,[Tu,Tv,x],Tt]+TD_\theta(Ts,[Tu,Tv,x])t-[Ts,x,T[u,v,t]_{T}]+TD_\theta(Ts,x)[u,v,t]_{T}\\
=&[Tu,Tv,[Ts,x,Tt]]-[Tu,Tv,TD_\theta(Ts,x)t]- [[Tu,Tv,Ts], x, Tt]+TD_\theta([Tu,Tv,Ts], x)t\\
 &-[Ts,[Tu,Tv,x],Tt]+TD_\theta(Ts,[Tu,Tv,x])t-[Ts,x,[Tu,Tv,Tt]]+TD_\theta(Ts,x)[u,v,t]_{T}\\
 =&T\big(-D_\theta(Tu,Tv)D_\theta(Ts,x)t+D_\theta([Tu,Tv,Ts], x)t+D_\theta(Ts,[Tu,Tv,x])t+D_\theta(Ts,x)D_\theta(Tu,Tv)t\big)\\
=&0,
%&\rho^l_E(u,v,\rho^r_E(x,s,t))-\rho^r_E(\rho^l_E(u,v,x), s,t)- \rho^r_E(x,[u,v,s]_{E},t)-\rho^r_E(x,s,[u,v,t]_{E})\\
%=&[Eu,Ev,[x,Es,Et]]-[Eu,Ev,ED_\theta(x,Es)t]-[[Eu,Ev,x], Es,Et]+ED_\theta([Eu,Ev,x], Es)t\\
%&- [x,E[u,v,s]_{E},Et]+ED_\theta(x,E[u,v,s]_{E})t-[x,Es,E[u,v,t]_{E}]+ED_\theta(x,Es)[u,v,t]_{E}\\
%=&E\big(-D_\theta(Eu,Ev)D_\theta(x,Es)t+D_\theta([Eu,Ev,x], Es)t+D_\theta(x,[Eu,Ev,Es])t+D_\theta(x,Es)D_\theta(Eu,Ev)t\big)\\
%=&0,
\end{align*}
which implies that Eqs. \eqref{4.2} holds. Similarly, we can deduce that Eqs. \eqref{4.3}-\eqref{4.5} hold. Therefore, we deduce that $(\mathbb{L};\rho^l_T,\rho^m_T,\rho^r_T)$  is a representation of the 3-Leibniz algebra $(\mathbb{L}',[-,-,-]_T)$.
\end{proof}

Let $T:\mathbb{L}'\rightarrow \mathbb{L}$ be   a nonabelian embedding tensor from a Lie triple system  $(\mathbb{L}',[-,-,-]')$  to another Lie triple system  $(\mathbb{L},[-,-,-])$ with respect to  a coherent action $\theta$.
Recall that Corollary \ref{coro:NET} and Lemma \ref{lemma:representation} give a new 3-Leibniz algebra  $(\mathbb{L}',[-,-,-]_T)$  and
a new representation   $(\mathbb{L};\rho^l_T,\rho^m_T,\rho^r_T)$  over    $(\mathbb{L}',[-,-,-]_T)$. Consider the cochain complex of  $(\mathbb{L}',[-,-,-]_T)$ with coefficients in  $(\mathbb{L};\rho^l_T,\rho^m_T,\rho^r_T)$ :
\begin{equation*}
(\mathcal{C}_{\mathrm{3Leib}}^{\bullet}(\mathbb{L}',\mathbb{L}),\partial_T^\bullet)=(\oplus_{n=1}^{\infty}\mathcal{C}_{\mathrm{3Leib}}^{n}(\mathbb{L}',\mathbb{L}),\partial_T^\bullet).
\end{equation*}
More precisely,
\begin{equation*}
\mathcal{C}_{\mathrm{3Leib}}^{n}(\mathbb{L}',\mathbb{L})=\mathrm{Hom}(\underbrace{\wedge^2 \mathbb{L}'\otimes\cdots\otimes\wedge^2 \mathbb{L}'}_{n-1} \otimes \mathbb{L}',\mathbb{L})
\end{equation*}
 and its coboundary operator $\partial_T^{n}:\mathcal{C}_{\mathrm{3Leib}}^{n}(\mathbb{L}',\mathbb{L})\rightarrow \mathcal{C}_{\mathrm{3Leib}}^{n+1}(\mathbb{L}',\mathbb{L}), f\mapsto  (\partial_T^{n}f)$ is given as follows:
\begin{align*}
&(\partial_T^{n} f)(\mathfrak{U}_1,\mathfrak{U}_2, \ldots, \mathfrak{U}_n,w)\\
=&\sum_{1\leq j<k\leq n}(-1)^jf(\mathfrak{U}_1,\ldots,\widehat{\mathfrak{U}_j},\ldots,\mathfrak{U}_{k-1},u_k \wedge[u_j,v_j,v_k]_{T}+[u_j,v_j,u_k]_{T}\wedge v_k,\ldots,\mathfrak{U}_n,w)\\
&+\sum_{j=1}^n(-1)^jf(\mathfrak{U}_1,\ldots,\widehat{\mathfrak{U}_j},\ldots,\mathfrak{U}_{n},[u_j,v_j,w]_{T})+\sum_{j=1}^n(-1)^{j+1}\rho^l_T(\mathfrak{U}_j,f(\mathfrak{U}_1,\ldots,\widehat{\mathfrak{U}_j},\ldots,\mathfrak{U}_{n},w))\\
&+(-1)^{n+1}(\rho^m_T(u_n, f(\mathfrak{U}_1,\cdots,\mathfrak{U}_{n-1},v_n),w)+\rho^r_T(f(\mathfrak{U}_1,\ldots,\mathfrak{U}_{n-1},u_n), v_n,w)),
\end{align*}
for $\mathfrak{U}_i=u_i\wedge v_i\in \wedge^2 \mathbb{L}', 1\leq i\leq n$ and $w\in \mathbb{L}'$.  In particular, for   $f\in \mathcal{C}_{\mathrm{3Leib}}^{1}(\mathbb{L}',\mathbb{L})$,
 \begin{align*}
&(\partial_T^{1} f)(u, v, w)\\
=& -f([u, v, w]_T)+\rho^l_T(u,v,f(w))+\rho^m_T(u,f(v),w)+\rho^r_T(f(u),v,w)\\
=&-f(D_\theta(Tu,Tv)w+[u, v, w]')+[Tu,Tv,f(w)]+[Tu,f(v),Tw]-TD_\theta(Tu,f(v))w\\
&+[f(u),Tv,Tw]-TD_\theta(f(u),Tv)w.
\end{align*}
 For  $(a,b)\in\mathcal{C}_{\mathrm{3Leib}}^{0}(\mathbb{L}',\mathbb{L}):=\wedge^2\mathbb{L}$, we define
 $\partial_T^{0}:\mathcal{C}_{\mathrm{3Leib}}^{0}(\mathbb{L}',\mathbb{L})\rightarrow \mathcal{C}_{\mathrm{3Leib}}^{1}(\mathbb{L}',\mathbb{L}), (a,b)\mapsto  \Im(a,b)$ by
\begin{align*}
 \Im(a,b)u=TD_\theta(a,b)u-[a,b,Tu], ~~\forall ~  u\in \mathbb{L}'.
\end{align*}

\begin{prop}\label{prop:1-cocycle}
  With the above notations,  we have  $\partial_T^{1}\Im(a,b)=0,$  that is the composition  $ \mathcal{C}_{\mathrm{3Leib}}^{0}(\mathbb{L}',\mathbb{L})\stackrel{\partial_T^{0}}{\longrightarrow} \mathcal{C}_{\mathrm{3Leib}}^{1}(\mathbb{L}',\mathbb{L})\stackrel{\partial_T^{1}}{\longrightarrow} \mathcal{C}_{\mathrm{3Leib}}^{2}(\mathbb{L}',\mathbb{L})$ is the zero map.
\end{prop}
\begin{proof}
For any $u,v,w\in V$,   by Eqs. \eqref{2.3}, \eqref{2.7}, \eqref{2.11},  \eqref{2.14} and \eqref{2.18}, we have
\begin{align*}
 & (\partial_T^{1}\Im(a,b))(u,v,w)\\
=&-\Im(a,b)D_\theta(Tu,Tv)w-\Im(a,b)[u, v, w]'+[Tu,Tv,\Im(a,b)w]+[Tu,\Im(a,b)v,Tw]-\\
&TD_\theta(Tu,\Im(a,b)v)w+[\Im(a,b)u,Tv,Tw]-ED_\theta(\Im(a,b)u,Tv)w\\
=&-T(D_\theta(a,b)D_\theta(Tu,Tv)w)+[a,b,TD_\theta(Tu,Tv)w]-TD_\theta(a,b)[u, v, w]'+[a,b,T[u, v, w]']+\\
&[Tu,Tv,TD_\theta(a,b)w]-[Tu,Tv,[a,b,Tw]]+[Tu,TD_\theta(a,b)v,Tw]-[Tu,[a,b,Tv],Tw]-\\
&TD_\theta(Tu,TD_\theta(a,b)v)w+TD_\theta(Tu,[a,b,Tv])w+[TD_\theta(a,b)u,Tv,Tw]-\\
&[[a,b,Tu],Tv,Tw]-TD_\theta(TD_\theta(a,b)u,Tv)w+TD_\theta([a,b,Tu],Tv)w\\
=&-T(D_\theta(a,b)D_\theta(Tu,Tv)w)+[Tu,Tv,TD_\theta(a,b)w]+[Tu,TD_\theta(a,b)v,Tw]-\\
&TD_\theta(Tu,TD_\theta(a,b)v)w+TD_\theta(Tu,[a,b,Tv])w+[TD_\theta(a,b)u,Tv,Tw]-\\
&TD_\theta(TD_\theta(a,b)u,Tv)w+TD_\theta([a,b,Tu],Tv)w
%=&-T(D_\theta(a,b)D_\theta(Tu,Tv)w)+T(D_\theta(Tu,Tv)D_\theta(a,b)w)+[Tu,TD_\theta(a,b)v,Tw]-TD_\theta(Tu,TD_\theta(a,b)v)w\\
%&+TD_\theta(Tu,[a,b,Tv])w+[TD_\theta(a,b)u,Tv,Tw]-TD_\theta(TD_\theta(a,b)u,Tv)w+TD_\theta([a,b,Tu],Tv)w\\
\end{align*}
\begin{align*}
=&T\Big(-D_\theta(a,b)D_\theta(Tu,Tv)w+D_\theta(Tu,Tv)D_\theta(a,b)w+[u,v,D_\theta(a,b)w]'+D_\theta(Tu,TD_\theta(a,b)v)w+\\
&[u,D_\theta(a,b)v, w]'-D_\theta(Tu,TD_\theta(a,b)v)w+D_\theta(Tu,[a,b,Tv])w+D_\theta(TD_\theta(a,b)u,Tv)w+\\
&[D_\theta(a,b)u,v,w]'-D_\theta(TD_\theta(a,b)u,Tv)w+D_\theta([a,b,Tu],Tv)w\Big)\\
=&0.
\end{align*}
 Thus, $\partial_T^{1}\Im(a,b)=0.$
\end{proof}

Next we define the cohomology theory of a    nonabelian embedding tensor $T$ from a Lie triple system  $(\mathbb{L}',[-,-,-]')$  to another Lie triple system  $(\mathbb{L},[-,-,-])$ with respect to  a coherent action $\theta$.
 For $n\geq 0$, define the set of $n$-cochains
of $T$ by
$$\mathcal{C}_{T}^{n}(\mathbb{L}',\mathbb{L})=\mathcal{C}_{\mathrm{3Leib}}^{n}(\mathbb{L}',\mathbb{L}).$$
Then, by Proposition \ref{prop:1-cocycle}, $(\oplus_{n=0}^{\infty}\mathcal{C}_{T}^{n}(\mathbb{L}',\mathbb{L}),\partial_T^\bullet)$ is a cochain complex. For $n\geq 1$, we
denote the set of $n$-cocycles by $\mathcal{Z}_{T}^{n}(\mathbb{L}',\mathbb{L})$, the set of $n$-coboundaries by $\mathcal{B}_{T}^{n}(\mathbb{L}',\mathbb{L})$ and the $n$-th
cohomology group of the  nonabelian  embedding tensor $T$ by $\mathcal{H}^n_{T}(\mathbb{L}',\mathbb{L})=\mathcal{Z}^n_{T}(\mathbb{L}',\mathbb{L})/\mathcal{B}^n_{T}(\mathbb{L}',\mathbb{L})$.

Furthermore, comparing the coboundary operator $\partial_T^{n}$ given above with the    differential $\mathrm{d}^{n}_T$
 given in Eq. \eqref{3.2}, we obtain the following result.
\begin{theorem}
 Let $T:\mathbb{L}'\rightarrow \mathbb{L}$ be   a nonabelian embedding tensor from a Lie triple system  $(\mathbb{L}',[-,-,-]')$  to another Lie triple system  $(\mathbb{L},[-,-,-])$ with respect to  a coherent action $\theta$.  Then  we have
$\partial_T^{n}f=(-1)^{n-1}\mathrm{d}^{n}_T(f)$ for $f\in \mathcal{C}_{T}^{n}(\mathbb{L}',\mathbb{L}), n=1,2,\ldots.$
\end{theorem}

\begin{proof}
For any $\mathfrak{U}_i=u_i\wedge v_i\in \wedge^2 V, 1\leq i\leq n$ and $w\in V$, according to Lemma \ref{lemma:ETLTS}, we have
\begin{align*}
 &\mathrm{d}_T^n(f)(\mathfrak{U}_1,\mathfrak{U}_2, \ldots, \mathfrak{U}_n,w)\\%=l_1^T(f)(\mathfrak{U}_1,\mathfrak{U}_2, \ldots, \mathfrak{U}_n,w)\\
 =&(l_1(f)+\frac{1}{2}l_3(T,T,f))(\mathfrak{U}_1,\mathfrak{U}_2, \ldots, \mathfrak{U}_n,w)\\
 =&[\varpi\boxtimes D_\theta\boxtimes \varpi',f]_{\mathrm{3LA}}(\mathfrak{U}_1,\mathfrak{U}_2, \ldots, \mathfrak{U}_n,w)+
 \frac{1}{2} [[[\varpi\boxtimes D_\theta\boxtimes \varpi',T]_{\mathrm{3LA}},T]_{\mathrm{3LA}},f]_{\mathrm{3LA}}(\mathfrak{U}_1,\mathfrak{U}_2, \ldots, \mathfrak{U}_n,w)\\
 %=&\frac{1}{2} \big([[\varpi\boxtimes D_\theta,E]_{\mathrm{NA}},E]_{\mathrm{NA}}\circ f-(-1)^{n-1}f\circ [[\varpi\boxtimes D_\theta,E]_{\mathrm{NA}},E]_{\mathrm{NA}}\big)(\mathfrak{U}_1,\mathfrak{U}_2, \ldots, \mathfrak{U}_n,w)\\
 =&-(-1)^{n-1}\sum_{j=1}^{n-1}\sum_{i=1}^j(-1)^{i+1}f(\mathfrak{U}_1,  \ldots,\widehat{\mathfrak{U}}_i, \ldots, \mathfrak{U}_j, [u_i, v_i,u_{j+1}]\wedge v_{j+1}, \mathfrak{U}_{j+2},\ldots, \mathfrak{U}_{n}, w)\\
 &-(-1)^{n-1}\sum_{j=1}^{n-1}\sum_{i=1}^j(-1)^{i+1}f(\mathfrak{U}_1,  \ldots,\widehat{\mathfrak{U}}_i, \ldots, \mathfrak{U}_j, u_{j+1}\wedge[u_i, v_i,v_{j+1}], \mathfrak{U}_{j+2},\ldots, \mathfrak{U}_{n}, w)\\
 &-(-1)^{n-1} \sum_{i=1}^n(-1)^{i+1}f(\mathfrak{U}_1,\ldots,\widehat{\mathfrak{U}}_i, \ldots, \mathfrak{U}_n,  [u_i, v_i,w]) \\
 &+(-1)^{n-1} \sum_{i=1}^n(-1)^{i+1}\rho_T^l(u_i,v_i,f(\mathfrak{U}_1,\ldots,\widehat{\mathfrak{U}}_i, \ldots, \mathfrak{U}_n,    w))
 \end{align*}
\begin{align*}
 &-(-1)^{n-1}\sum_{j=1}^{n-1}\sum_{i=1}^j(-1)^{i+1}f(\mathfrak{U}_1,  \ldots,\widehat{\mathfrak{U}}_i, \ldots, \mathfrak{U}_j, D_\theta(Tu_i,Tv_i)u_{j+1}\wedge v_{j+1}, \mathfrak{U}_{j+2},\ldots, \mathfrak{U}_{n}, w)\\
 &-(-1)^{n-1}\sum_{j=1}^{n-1}\sum_{i=1}^j(-1)^{i+1}f(\mathfrak{U}_1,  \ldots,\widehat{\mathfrak{U}}_i, \ldots, \mathfrak{U}_j, u_{j+1}\wedge D_\theta(Tu_i,Tv_i)v_{j+1}, \mathfrak{U}_{j+2},\ldots, \mathfrak{U}_{n}, w)\\
 &-(-1)^{n-1} \sum_{i=1}^n(-1)^{i+1}f(\mathfrak{U}_1,\ldots,\widehat{\mathfrak{U}}_i, \ldots, \mathfrak{U}_n,D_\theta(Tu_i,Tv_i)w) \\
 &+\rho_T^m(u_n,f(\mathfrak{U}_1,\ldots,\widehat{\mathfrak{U}}_i, \ldots, \mathfrak{U}_{n-1},  v_n),w)+\rho_T^r(f(\mathfrak{U}_1,\ldots,\widehat{\mathfrak{U}}_i, \ldots, \mathfrak{U}_{n-1},  u_n),v_n,w)\\
=&(-1)^{n-1}\partial_T^{n}f,
\end{align*}
which implies that $\partial_T^{n}f=(-1)^{n-1}\mathrm{d}^{n}_T(f)$.
The proof is completed.
\end{proof}

\section{  Deformations of nonabelian embedding tensors between Lie triple systems }\label{sec: Deformations}
\def\theequation{\arabic{section}.\arabic{equation}}
\setcounter{equation} {0}

In this section, we study infinitesimal deformations of nonabelian embedding tensors between Lie triple systems using the cohomology theory given in the previous section, following Gerstenhaber's approach.
In particular, we introduce the notion of a Nijenhuis element associated with a nonabelian embedding tensor%  that arises from trivial infinitesimal deformations
, which
leads to a trivial infinitesimal deformations.

Let $T:\mathbb{L}'\rightarrow \mathbb{L}$ be   a nonabelian embedding tensor from a Lie triple system  $(\mathbb{L}',[-,-,-]')$  to another Lie triple system  $(\mathbb{L},[-,-,-])$ with respect to  a coherent action $\theta$.

 \begin{definition}
A parameterized sum $T_t=T+tT_1$, for some $T_1\in \mathrm{Hom}(\mathbb{L}',\mathbb{L})$, is called an infinitesimal deformation of $T$ if $T_t$ is a nonabelian embedding tensor
 for all values of parameter $t$. In this case, we say that $T_1$ generates an infinitesimal deformation of $T$.
\end{definition}

Suppose that $T_1$  generates an infinitesimal deformation of  $T$, then we have
$$ [T_tu,T_tv,T_tw]=T_t(D_\theta(T_tu,T_tv)w+[u,v,w]'),$$
for all $u,v,w\in \mathbb{L}'$. This is equivalent to the following equations
 \begin{align}
&[Tu,Tv,T_1w]+[Tu,T_1v,Tw]+[T_1u,Tv,Tw]\nonumber\\
&=T_1(D_\theta(Tu,Tv)w+[u,v,w]')+T(D_\theta(T_1u,Tv)w+D_\theta(Tu,T_1v)w),\label{5.1}\\
&[Tu,T_1v,T_1w]+[T_1u,Tv,T_1w]+[T_1u,T_1v,Tw]\nonumber\\
&=TD_\theta(T_1u,T_1v)w+T_1(D_\theta(T_1u,Tv)w+D_\theta(Tu,T_1v)w),\label{5.2}\\
&[T_1u,T_1v,T_1w]=T_1D_\theta(T_1u,T_1v)w.\label{5.3}
\end{align}
Thus, $T_t$ is a infinitesimal deformation of $T$ if and only if Eqs. \eqref{5.1}-\eqref{5.3} hold.
Observe that  Eq. \eqref{5.1} implies that $T_1$ is a 1-cocycle the nonabelian embedding tensor $T$.
From Eq. \eqref{5.3} it follows that the map $T_1$ is an embedding tensor  on the
Lie triple system  $(\mathbb{L},[-,-,-])$ with respect to the representation $(\mathbb{L}',\theta)$ (see \cite{Teng25}).

Let  $(\mathfrak{g},[-,-,-]_{\mathfrak{g}})$  be a 3-Leibniz algebra and  $\omega$ be a linear map. If for any values of parameter $t$, the
multiplication $[-,-,-]_t$ defined by
$$[u,v,w]_t=[u,v,w]+t\omega(u,v,w)$$
for $u,v,w\in \mathfrak{g}$,  also gives a 3-Leibniz algebra structure, we say that $\omega$ generates an infinitesimal deformation of
the 3-Leibniz algebra $(\mathfrak{g},[-,-,-]_{\mathfrak{g}})$ . The two types of infinitesimal deformations are related as follows.

  \begin{prop}
If $T_1$  generates an infinitesimal deformation of  $T$, then the product $\omega_1$ on $\mathbb{L}'$ defined by
$$\omega_1(u,v,w)=D_\theta(T_1u,Tv)w+D_\theta(Tu,T_1v)w$$
for $u,v,w\in  \mathbb{L}'$, generates an infinitesimal deformation of the descendent 3-Leibniz algebra   $(\mathbb{L}',[-,-,-]_T)$ given  in
Corollary \ref{coro:NET}.
\end{prop}

\begin{proof}
According to Corollary \ref{coro:NET}, we denote by $[-,-,-]_t$
the descendent  3-Leibniz algebra structure corresponding  to the nonabelian embedding tensor $T_t$.
Then, for any $u,v,w\in  \mathbb{L}'$, we have
  \begin{align*}
&[u,v,w]_t= D_\theta(T_tu,T_tv)w+[u,v,w]'\\
&=D_\theta(Tu,Tv)w+[u,v,w]'+t(D_\theta(T_1u,Tv)w+D_\theta(Tu,T_1v)w)~~~( \mathrm{mod}~t^2)\\
&=[u,v,w]_T+t\omega_1(u,v,w),
\end{align*}
which implies that $\omega_1$ generates an infinitesimal deformation of  $(\mathbb{L}',[-,-,-]_T)$.
\end{proof}

 \begin{definition}
Two infinitesimal deformations $T_t=T+tT_1$ and $\widetilde{T}_t=T+t\widetilde{T}_1$
of a  nonabelian embedding tensor $T$ from a Lie triple system  $(\mathbb{L}',[-,-,-]')$  to another Lie triple system  $(\mathbb{L},[-,-,-])$ with respect to  a coherent action $\theta$  are said to
be equivalent if there exists an element $\mathfrak{A}=(a, b)\in \wedge^2\mathbb{L}$ such that the pair
$$(\phi_t=\mathrm{id}_{\mathbb{L}}+t[a,b,-],\phi'_t=\mathrm{id}_{\mathbb{L}'}+tD_\theta(a,b))$$
defines a homomorphism of nonabelian embedding tensors
from $\widetilde{T}_t$ to $T_t$.
In particular, an infinitesimal deformation $T_t=T+tT_1$ of a nonabelian embedding
tensor $T$ is called trivial if there exists an element  $(a, b)\in \wedge^2\mathbb{L}$ such that $(\phi_t,\phi'_t)$ is a
homomorphism from $T_t$
to $T$.
 \end{definition}

Let $(\phi_t,\phi'_t)$ be a homomorphism from from $\widetilde{T}_t$ to $T_t$.  Then $\phi_t$ and $\phi'_t$ are   Lie triple system  homomorphisms, which implies that
\begin{small}
\begin{align}
\left\{ \begin{array}{lll}
~[x,[a,b,y],[a,b,z]]+[[a,b,x],y,[a,b,z]]+[[a,b,x],[a,b,y],z]=0,\\
 ~[[a,b,x],[a,b,y],[a,b,z]]=0,\\
 ~[u,D_\theta(a,b)v,D_\theta(a,b)w]'+[D_\theta(a,b)u,v,D_\theta(a,b)w]'+[D_\theta(a,b)u,D_\theta(a,b)v,w]'=0,\\
 ~[D_\theta(a,b)u,D_\theta(a,b)v,D_\theta(a,b)w]'=0,
 \end{array}  \right. \label{5.4}
 \end{align}
 \end{small}
for all $x,y,z\in \mathbb{L}$ and $u,v,w\in \mathbb{L}'$.

Note that by Eqs. \eqref{2.16} and \eqref{2.17}, we obtain that for all $x,y\in \mathbb{L}$ and $u\in \mathbb{L}'$,
\begin{align}
\left\{ \begin{array}{lll}
\theta([a,b,x],[a,b,y])u+\theta([a,b,x],y)D_\theta(a,b)u+\theta(x,[a,b,y])D_\theta(a,b)u=0,\\
 \theta([a,b,x],[a,b,y])D_\theta(a,b)u=0,\\
 D_\theta([a,b,x],[a,b,y])u+D_\theta([a,b,x],y)D_\theta(a,b)u+D_\theta(x,[a,b,y])D_\theta(a,b)u=0,\\
D_\theta([a,b,x],[a,b,y])D_\theta(a,b)u=0.
 \end{array}  \right. \label{5.5}
 \end{align}

Moreover, Eq.   \eqref{2.15}   yields that for all  $u\in \mathbb{L}'$,
 \begin{align*}
&\widetilde{T}_tu+t[a,b,\widetilde{T}_tu]=T_t(u+tD_\theta(a,b)u),
\end{align*}
which implies
  \begin{align}
\widetilde{T}_1u+[a,b,Tu]&=T_1u +T(D_\theta(a,b)u), \label{5.6}\\
[a,b,\widetilde{T}_1u]&=T_1D_\theta(a,b)u.\nonumber
\end{align}
Note that Eq. \eqref{5.6}  means that $\widetilde{T}_1u-T_1u= \partial_T^0(a,b)u$.
Thus, we have the following theorem.

\begin{theorem}
%Let $T:\mathbb{L}'\rightarrow \mathbb{L}$ be   a nonabelian embedding tensor from a Lie triple system  $(\mathbb{L}',[-,-,-]')$  to another Lie triple system  $(\mathbb{L},[-,-,-])$ with respect to  a coherent action $\theta$.
 If $T_t=T+tT_1$ and $\widetilde{T}_t=T+t\widetilde{T}_1$
are two equivalent infinitesimal deformations
of a   nonabelian embedding tensor $T$, then $T_1$ and $\widetilde{T}_1$
define the same cohomology class in  $\mathcal{H}_T^1(\mathbb{L}',\mathbb{L})$.
\end{theorem}

 \begin{definition}
Let $T:\mathbb{L}'\rightarrow \mathbb{L}$ be   a nonabelian embedding tensor from a Lie triple system  $(\mathbb{L}',[-,-,-]')$  to another Lie triple system  $(\mathbb{L},[-,-,-])$ with respect to  a coherent action $\theta$.
 An element $\mathfrak{A}=(a, b)\in \wedge^2\mathbb{L}$ is said to be a Nijenhuis element associated to $T$ if $\mathfrak{A}$
satisfies Eqs. \eqref{5.4} and \eqref{5.5}  and the equation
  \begin{align}
&[a,b,T(D_\theta(a,b)u)-[a,b,Tu]]=0 ,\label{5.7}
\end{align}
for all $u\in  \mathbb{L}'$.  Denote by $\mathrm{Nij}(T)$ the set of Nijenhuis elements associated with a nonabelian embedding tensor $T$.
 \end{definition}

It is easy to see that a trivial infinitesimal deformation of a nonabelian embedding tensor  between  Lie triple systems gives rise to a Nijenhuis element. However, the converse is also true.

 \begin{theorem}
Let $T:\mathbb{L}'\rightarrow \mathbb{L}$ be   a nonabelian embedding tensor from a Lie triple system  $(\mathbb{L}',[-,-,-]')$  to another Lie triple system  $(\mathbb{L},[-,-,-])$ with respect to  a coherent action $\theta$. Then, for any $\mathfrak{A}=(a, b)\in \mathrm{Nij}(T)$, $T_t=T+tT_1$ with $T_1=\partial_T^0(a,b)$ is
a trivial linear deformation of the nonabelian embedding tensor $T$.
\end{theorem}

 \begin{proof}
 For any Nijenhuis element $\mathfrak{A}=(a, b)\in \mathrm{Nij}(T)$, $T_t=T+tT_1$, we define a linear map $T_1:\mathbb{L}'\rightarrow\mathbb{L}'$ by
$$T_1u=\partial_T^0(a,b)u= T(D_\theta(a,b)u)-[a,b,Tu],$$
for all $u\in \mathbb{L}'$.
Let $T_t=T+tT_1$. By the definition of Nijenhuis elements, for all $x,y,z\in \mathbb{L}$ and  $u,v,w\in \mathbb{L}'$, we have
 \begin{align*}
[x,y,z]+t[a,b,[x,y,z]]&=[x+t[a,b,x],y+t[a,b,y],z+t[a,b,z]],\\
[u,v,w]'+tD_\theta(a,b)[u,v,w]'&=[u+tD_\theta(a,b)u,v+tD_\theta(a,b)v,w+tD_\theta(a,b)w]',\\
T_tu+t[a,b,T_tu]&=Tu+tTD_\theta(a,b)u,\\
D_\theta(x,y)u+tD_\theta(a,b)D_\theta(x,y)u&=D_\theta(x+[a,b,x],y+[a,b,y])(u+tD_\theta(a,b)u).
\end{align*}
Since $[a,b,-]$ and $D_\theta(a,b)$ are linear transformations of finite-dimensional  vector spaces $\mathbb{L}$ and $\mathbb{L}'$, respectively, it follows that
$\phi_t=\mathrm{id}_{\mathbb{L}}+t[a,b,-]$ and $\phi'_t=\mathrm{id}_{\mathbb{L}'}+tD_\theta(a,b)$ are Lie triple system isomorphisms for $|t|$  sufficiently small.
Thus, by  Proposition \ref{prop:newnonabelian},
 \begin{align*}
& \phi_t^{-1}\circ T\circ\phi'_t=(\mathrm{id}_{\mathbb{L}}+t[a,b,-])^{-1}\circ T\circ(\mathrm{id}_{\mathbb{L}'}+tD_\theta(a,b))\\
&=\sum_{i=0}^{+\infty}(-t[a,b,-])^i\circ(T+tTD_\theta(a,b))\\
&=T+t(TD_\theta(a,b)-[a,b,T-])+\sum_{i=1}^{+\infty}(-1)^it^{i+1}[a,b,-]^i\circ(TD_\theta(a,b)-[a,b,T-])\\
&=T+t(TD_\theta(a,b)-[a,b,T-])\\
&=T+tT_1=T_t
\end{align*}
is a nonabelian embedding tensor from a Lie triple system  $(\mathbb{L}',[-,-,-]')$  to another Lie triple system  $(\mathbb{L},[-,-,-])$ with respect to  a coherent action $\theta$ for $|t|$  sufficiently small.
Thus, $T_1$   satisfies the Eqs.  \eqref{5.1}-\eqref{5.3}. Therefore, $T_t$
is a nonabelian embedding tensor for all $t$, which means that $T_1$  generates
an infinitesimal deformation of $T$. It is straightforward to show that this infinitesimal deformation is
trivial.
\end{proof}

 We end this section by presenting a concrete example of an infinitesimal  deformation of a
nonabelian embedding tensor $T$ on the  3-dimensional Lie triple system  $(\mathbb{L}_3,[-,-,-])$  given in Example \ref{exam:3-dimensional}.

 \begin{exam}
 Let $(\mathbb{L}_3,[-,-,-])$ be the  3-dimensional Lie triple system  given in Example \ref{exam:3-dimensional}. Then,
 $$T_t=\begin{bmatrix}
 0 &  \frac{1}{a} & 0\\
 a & c & 0\\
 b & d & t
 \end{bmatrix}$$
   is an infinitesimal  deformation of the nonabelian embedding
tensor $T=\begin{bmatrix}
  0 &  \frac{1}{a} & 0\\
 a & c & 0\\
 b & d & 0
 \end{bmatrix}$ on $(\mathbb{L}_3,[-,-,-])$ with respect to the adjoint action  for all $a\neq0,b,c,d\in \mathbf{k}$.
\end{exam}

\section{ From   nonabelian embedding tensors between Lie algebras to those between Lie
triple systems }\label{sec: From}
\def\theequation{\arabic{section}.\arabic{equation}}
\setcounter{equation} {0}

Motivated by the construction of Lie triple systems from Lie algebras. We give some connections between nonabelian embedding tensors on Lie algebras and Lie triple systems.

\begin{definition} \cite{Tang}
Let $(L,[-,-])$ and $(L',[-,-]')$ be two Lie algebras. A Lie algebra homomorphism
$\rho:L\rightarrow \mathrm{Der}_{LA}(L')$ is called an action of $L$ on $L'$. An action $\rho$ of $L$ on $L'$ is called a
coherent action if it satisfies
\begin{align}
&  [\rho(x)u,v]'=0,\label{6.1}
\end{align}
for all $x\in L$ and $u,v\in L'$.
\end{definition}

\begin{definition} \cite{Tang}
A nonabelian embedding tensor from a Lie algebra $(L',[-,-]')$ to another Lie algebra $(L,[-,-])$  with respect to a
coherent action $\rho$ is a linear map $T:L'\rightarrow L$ such that
\begin{align}
&  [Tu,Tv]=T(\rho(Tu)v+[u,v]'),\label{6.2}
\end{align}
for all $ u, v\in L'$.
\end{definition}

As we know, there is a method of constructing Lie
triple systems from Lie algebras given in Remark \ref{remark: Lie-LTS}. In the
sequel, we study the relationship between Lie algebras and the corresponding Lie
triple systems.

  \begin{prop}
 Let $\rho$ be a coherent action of a Lie algebra $(L,[-,-])$ on
another Lie algebra $(L',[-,-]')$. Define $\theta_\rho:L\times L\rightarrow \mathrm{End}(L')$ by
\begin{align}
& \theta_\rho(x,y)u=\rho(y)\rho(x)u,\label{6.3}
\end{align}
for all $x,y\in L$ and $u\in L'$. Then $\theta_\rho$ is a coherent action of a Lie triple system $(L,[-,-,-]=[-,-]\circ([-,-]\otimes \mathrm{id}))$ on
another Lie triple system  $(L',[-,-,-]'=[-,-]'\circ([-,-]'\otimes \mathrm{id}))$.
\end{prop}

\begin{proof}
According to \cite{Hajjaji24}, we get that $\theta_\rho$ is a  action of a Lie triple system $(L,[-,-,-])$ on
another Lie triple system  $(L',[-,-,-]')$.
Moreover, for any $x,y\in L$, using $\rho(x)\in \mathrm{Der}_{LA}(L')$ and Eq. \eqref{6.1}, we have
%\begin{align*}
%& \theta_\rho(x,y)[u,v,w]'=\rho(y)\rho(x) [[u,v]',w]'\\
%&=\rho(y)([\rho(x)[u,v]',w]'+[[u,v]',\rho(x)w]')\\
%&=\rho(y)([[\rho(x)u,v]',w]'+[[u,\rho(x)v]',w]'+[[u,v]',\rho(x)w]')\\
%&=[[\rho(y)\rho(x)u,v]',w]'+[[\rho(x)u,\rho(y)v]',w]'+[[\rho(x)u,v]',\rho(y)w]'+[[\rho(y)u,\rho(x)v]',w]'+\\
%&~~~~[[u,\rho(y)\rho(x)v]',w]'+[[[u,\rho(x)v]',\rho(y)w]'+[[\rho(y)u,v]',\rho(x)w]'+[[u,\rho(y)v]',\rho(x)w]'+\\
%&~~~~[[u,v]',\rho(y)\rho(x)w]'\\
%&=[\rho(y)\rho(x)u,v,w]'+[u,\rho(y)\rho(x)v,w]'+[u,v,\rho(y)\rho(x)w]'\\
%&=[\theta_\rho(x,y)u,v,w]'+[u,\theta_\rho(x,y)v,w]'+[u,v,\theta_\rho(x,y)w]'.
%\end{align*}
 $\theta_\rho(x,y)\in \mathrm{Der}_{Lts}(L')$.
\end{proof}

  \begin{prop}
 Let $T:L'\rightarrow L$  be a  nonabelian embedding tensor from a Lie algebra $(L',[-,-]')$ to another Lie algebra $(L,[-,-])$  with respect to a
coherent action $\rho$.
Then $T$ is a  nonabelian embedding tensor from a  Lie triple system $(L',[-,-,-]'=[-,-]'\circ([-,-]'\otimes \mathrm{id}))$ to
another Lie triple system $(L,[-,-,-]=[-,-]\circ([-,-]\otimes \mathrm{id}))$. with respect to a
coherent action $\theta_\rho$.
\end{prop}

\begin{proof}
For any $u,v,w\in L'$, by  Eqs. \eqref{2.6},\eqref{6.1}-\eqref{6.3}, we have
\begin{align*}
&[Tu,Tv,Tw]=[[Tu,Tv],Tw]=[T(\rho(Tu)v+[u,v]'),Tw]\\
&=T\big(\rho(T(\rho(Tu)v+[u,v]'))w+[\rho(Tu)v+[u,v]',w]'\big)\\
&=T\big(\rho(T(\rho(Tu)v))w+\rho(T[u,v]')w+[\rho(Tu)v,w]'+[[u,v]',w]'\big)\\
&=T\big(\rho([Tu,Tv]-T[u,v]')w+\rho(T[u,v]')w+[[u,v]',w]'\big)\\
&=T\big(\rho(Tu)\rho(Tv)w-\rho(Tv)\rho(Tu)w+[[u,v]',w]'\big)\\
&=T\big(\theta_\rho(Tv,Tu)w-\theta_\rho(Tu,Tv)w+[[u,v]',w]'\big)\\
&=T\big(D_{\theta_\rho}(Tu,Tv)w+[u,v,w]'\big),
\end{align*}
which implies that $T$  is a nonabelian embedding tensor  from  a  Lie triple system $(L',[-,-,-]')$ to
another Lie triple system $(L,[-,-,-])$. with respect to a coherent action $\theta_\rho$.
\end{proof}

%\begin{center}
% {\bf Acknowledgment}
% \end{center}
%
% The paper is supported by the Foundation of Science and Technology of Guizhou Province (Grant No. [2018]1020).

%\bibliographystyle{amsplain}

\end{document}